\documentclass[12pt]{article}
\usepackage{latexsym, amssymb}

\DeclareMathAlphabet\gothic{U}{euf}{m}{n}

\makeatletter
\def\eqnarray{\stepcounter{equation}\let\@currentlabel=\theequation
\global\@eqnswtrue
\tabskip\@centering\let\\=\@eqncr
$$\halign to \displaywidth\bgroup\hfil\global\@eqcnt\z@
  $\displaystyle\tabskip\z@{##}$&\global\@eqcnt\@ne
  \hfil$\displaystyle{{}##{}}$\hfil
  &\global\@eqcnt\tw@ $\displaystyle{##}$\hfil
  \tabskip\@centering&\llap{##}\tabskip\z@\cr}

\def\endeqnarray{\@@eqncr\egroup
      \global\advance\c@equation\m@ne$$\global\@ignoretrue}

\def\@yeqncr{\@ifnextchar [{\@xeqncr}{\@xeqncr[5pt]}}
\makeatother

\begin{document}
\bibliographystyle{tom}

\newtheorem{lemma}{Lemma}[section]
\newtheorem{thm}[lemma]{Theorem}
\newtheorem{cor}[lemma]{Corollary}
\newtheorem{voorb}[lemma]{Example}
\newtheorem{rem}[lemma]{Remark}
\newtheorem{prop}[lemma]{Proposition}
\newtheorem{stat}[lemma]{{\hspace{-5pt}}}

\newenvironment{remarkn}{\begin{rem} \rm}{\end{rem}}
\newenvironment{exam}{\begin{voorb} \rm}{\end{voorb}}

\newcommand{\gota}{\gothic{a}}
\newcommand{\gotb}{\gothic{b}}
\newcommand{\gotc}{\gothic{c}}
\newcommand{\gote}{\gothic{e}}
\newcommand{\gotf}{\gothic{f}}
\newcommand{\gotg}{\gothic{g}}
\newcommand{\gothh}{\gothic{h}}
\newcommand{\gotk}{\gothic{k}}
\newcommand{\gotm}{\gothic{m}}
\newcommand{\gotn}{\gothic{n}}
\newcommand{\gotp}{\gothic{p}}
\newcommand{\gotq}{\gothic{q}}
\newcommand{\gotr}{\gothic{r}}
\newcommand{\gots}{\gothic{s}}
\newcommand{\gotu}{\gothic{u}}
\newcommand{\gotv}{\gothic{v}}
\newcommand{\gotw}{\gothic{w}}
\newcommand{\gotz}{\gothic{z}}
\newcommand{\gotA}{\gothic{A}}
\newcommand{\gotB}{\gothic{B}}
\newcommand{\gotG}{\gothic{G}}
\newcommand{\gotL}{\gothic{L}}
\newcommand{\gotS}{\gothic{S}}
\newcommand{\gotT}{\gothic{T}}

\newcounter{teller}
\renewcommand{\theteller}{\Roman{teller}}
\newenvironment{tabel}{\begin{list}%
{\rm \bf \Roman{teller}.\hfill}{\usecounter{teller} \leftmargin=1.1cm
\labelwidth=1.1cm \labelsep=0cm \parsep=0cm}
                      }{\end{list}}

\newcounter{tellerr}
\renewcommand{\thetellerr}{\roman{tellerr}}
\newenvironment{subtabel}{\begin{list}%
{\rm  \roman{tellerr}.\hfill}{\usecounter{tellerr} \leftmargin=1.1cm
\labelwidth=1.1cm \labelsep=0cm \parsep=0cm}
                         }{\end{list}}

\newcounter{proofstep}
\newcommand{\nextstep}{\refstepcounter{proofstep}\ruimte \par 
          \noindent{\bf Step \theproofstep} \hspace{5pt}}
\newcommand{\firststep}{\setcounter{proofstep}{0}\nextstep}

\newcommand{\Ni}{{\bf N}}
\newcommand{\Ri}{{\bf R}}
\newcommand{\Ci}{{\bf C}}
\newcommand{\Ti}{{\bf T}}
\newcommand{\Zi}{{\bf Z}}
\newcommand{\Fi}{{\bf F}}

\newcommand{\proof}{\mbox{\bf Proof} \hspace{5pt}} 
\newcommand{\remark}{\mbox{\bf Remark} \hspace{5pt}}
\newcommand{\ruimte}{\vskip10.0pt plus 4.0pt minus 6.0pt}

\newcommand{\ad}{{\mathop{\rm ad}}}
\newcommand{\Ad}{{\mathop{\rm Ad}}}
\newcommand{\Aut}{\mathop{\rm Aut}}
\newcommand{\arccot}{\mathop{\rm arccot}}
\newcommand{\diam}{\mathop{\rm diam}}
\newcommand{\divv}{\mathop{\rm div}}
\newcommand{\codim}{\mathop{\rm codim}}
\newcommand{\RRe}{\mathop{\rm Re}}
\newcommand{\IIm}{\mathop{\rm Im}}
\newcommand{\Tr}{{\mathop{\rm Tr}}}
\newcommand{\Vol}{{\mathop{\rm Vol}}}
\newcommand{\card}{{\mathop{\rm card}}}
\newcommand{\supp}{\mathop{\rm supp}}
\newcommand{\sgn}{\mathop{\rm sgn}}
\newcommand{\essinf}{\mathop{\rm ess\,inf}}
\newcommand{\esssup}{\mathop{\rm ess\,sup}}
\newcommand{\Int}{\mathop{\rm Int}}
\newcommand{\Leibniz}{\mathop{\rm Leibniz}}
\newcommand{\lcm}{\mathop{\rm lcm}}
\newcommand{\loc}{{\rm loc}}
\newcommand{\rlim}{\mathop{\rm r.lim}}

\newcommand{\mod}{\mathop{\rm mod}}
\newcommand{\spann}{\mathop{\rm span}}
\newcommand{\ubar}{\underline{\;}}
\newcommand{\one}{1\hspace{-4.5pt}1}

\hyphenation{groups}
\hyphenation{unitary}

\newcommand{\cb}{{\cal B}}
\newcommand{\cc}{{\cal C}}
\newcommand{\cd}{{\cal D}}
\newcommand{\ce}{{\cal E}}
\newcommand{\cf}{{\cal F}}
\newcommand{\ch}{{\cal H}}
\newcommand{\ci}{{\cal I}}
\newcommand{\ck}{{\cal K}}
\newcommand{\cl}{{\cal L}}
\newcommand{\cm}{{\cal M}}
\newcommand{\co}{{\cal O}}
\newcommand{\cs}{{\cal S}}
\newcommand{\ct}{{\cal T}}
\newcommand{\cx}{{\cal X}}
\newcommand{\cy}{{\cal Y}}
\newcommand{\cz}{{\cal Z}}

\newfont{\fontcmrten}{cmr10}
\newcommand{\slbrl}{\mbox{\fontcmrten (}}
\newcommand{\slbrr}{\mbox{\fontcmrten )}}

 \thispagestyle{empty}

\begin{center}
{\Large{\bf Uniform subellipticity}}\\[2mm] 

\large  A.F.M. ter Elst$^1$ and  Derek W. Robinson$^2$\\[2mm]
\end{center}

\vspace{15mm}

\begin{center}
{\bf Abstract}
\end{center}

\begin{list}{}{\leftmargin=1.8cm \rightmargin=1.8cm \listparindent=10mm 
   \parsep=0pt}
\item
We establish two global  subellipticity properties of
positive symmetric second-order partial differential operators 
on $L_2(\Ri^d)$.
First, if $m \in \Ni$ then we consider operators $H_0$ with  coefficients in $W^{m+1,\infty}(\Ri^d)$ 
and domain $D(H_0)=W^{\infty,2}(\Ri^d)$ satisfying the subellipticity property
\[
c\,(\varphi, (I+H_0)\varphi)\geq \|\Delta^{\gamma/2} \varphi\|_2^2
\]
for some $c>0$ and $\gamma\in\langle0,1]$, uniformly for all $\varphi\in W^{\infty,2}(\Ri^d)$,
where $\Delta$ denotes the usual Laplacian.
Then we prove that $D(H^\alpha) \subseteq D(\Delta^{\alpha \gamma})$ for all $\alpha \in [0,2^{-1} (m + 1 + \gamma^{-1}) \rangle$.
Hence
there is a $c>0$ such that 
the  norm estimate
\[
c \, \|(I+H)^\alpha \varphi\|_2\geq \|\Delta^{\alpha \gamma} \varphi\|_2
\]
is valid for all $\varphi\in D(H^\alpha)$ where $H$ denotes the self-adjoint closure of $H_0$.
In particular, if the coefficients of $H_0$ are in $C_b^\infty(\Ri^d)$ then the conclusion is valid
for all $\alpha\geq0$.

Secondly, we prove that if 
\[
H_0=\sum^N_{i=1}X_i^* \, X_i
\;\;\; ,
\]
where the $X_i$ are vector fields on $\Ri^d$ with coefficients in $C_b^\infty(\Ri^d)$
satisfying a uniform version of 
H\"ormander's criterion for hypoellipticity, then $H_0$ satisfies the subellipticity 
condition for $\gamma=r^{-1}$ where $r$ is the rank of the set of vector fields.
Consequently 
$D(H^n) \subseteq D(\Delta^{n/r})$ for all $n \in \Ni$,
where $H$ is the closure of $H_0$.
\end{list}

\vfill

\noindent
October 2006

\vspace{0.5cm}

\noindent
AMS Subject Classification: 47B47, 47B44, 58G03.

\vspace{1cm}

\noindent
\begin{tabular}{@{}cl@{\hspace{10mm}}cl}
1. &Department of Mathematics   & 
  2. &Centre for Mathematics and its Applications  \\
& University of Auckland  & 
  &   Mathematical Sciences Institute  \\
&  Private bag 92019 & 
  & Australian National University \\
& Auckland & 
  & Canberra, ACT 0200\\
  & New Zealand.  & 
  & Australia. 
\end{tabular}

\newpage
\setcounter{page}{1}
\section{Introduction}\label{Shorf1}

Our aim is to derive two global  subellipticity properties of
second-order self-adjoint elliptic operators on $L_2(\Ri^d)$.
Initially we consider operators of the form
\begin{equation}
H_0=-\sum^d_{i,j=0}\partial_i\,c_{ij}\,\partial_j
\label{esub1.0}
\end{equation}
with domain $D(H_0) = W^{\infty,2}(\Ri^d)$,
where $\partial_0=iI$ and  $\partial_j=\partial/\partial x_j$ if $j\in\{1,\ldots,d\}$.
We  assume throughout that the coefficients $c_{ij}\in W^{m+1,\infty}(\Ri^d)$, where $m \in \Ni$, are complex-valued
and  $C=(c_{ij})$ is a symmetric positive-definite matrix.
In particular, the coefficients are always at least twice differentiable.
Although we allow the $c_{ij}$ to be complex one could use symmetry to re-express $H_0$ in the 
form (\ref{esub1.0}) but with real-valued coefficients. 
Then, however, the corresponding $c_{i0}$ and $c_{0j}$ are not necessarily in $W^{m+1,\infty}(\Ri^d)$.
Since $c_{ij}\in W^{2,\infty}(\Ri^d)$ it follows, however, that
$H_0$ is essentially self-adjoint on
$W^{\infty,2}(\Ri^d)$ (see, for example, \cite{Rob7},  Section~6, or Proposition~\ref{pcom15} below)
and we denote the self-adjoint closure by $H$.

If  $\gamma \in \langle0,1]$ then $H_0$ is defined to be {\bf subelliptic of order $\gamma$}
if  there is a $c>0$ such that
\begin{equation}
c\,(\varphi, (I+H_0)\varphi)\geq \|\Delta^{\gamma/2} \varphi\|_2^2
\label{esub1.1}
\end{equation}
for all $\varphi\in W^{\infty,2}(\Ri^d)$.
Then the subellipticity condition extends to $H$ and 
$c\,(I+H) \geq \Delta^\gamma$ in the sense of quadratic forms.
A local version of Condition~(\ref{esub1.1}) arose in H\"ormander's work \cite{Hor1} and is significant as it implies
hypoellipticity of $H_0$.
The global version implies  uniform boundedness of the semigroup kernel associated with 
$H$ by an argument based on Nash inequalities.

Our first result establishes that the subellipticity condition is self-improving.

\begin{thm}\label{tcom31}
Let $H_0$ be a positive, symmetric, subelliptic operator of order $\gamma \in \langle0,1]$
with coefficients $c_{ij} \in W^{m+1,\infty}(\Ri^d)$, where  $m \in \Ni$,
and with self-adjoint closure $H$.
Then $D(H^\alpha) \subseteq D(\Delta^{\alpha \gamma})$
for all  $\alpha \in [0,2^{-1} (m + 1 + \gamma^{-1}) \rangle$ and there is  a $c>0$ such that
\begin{equation}
c \, \|(I+H)^\alpha \varphi\|_2\geq \|\Delta^{\alpha \gamma} \varphi\|_2
\label{esub1.11}
\end{equation}
for all $\varphi\in D(H^\alpha)$.
\end{thm}

The theorem is a strengthened global version of a local result
of Fefferman and Phong (see \cite{FP}, first part of Theorem~1).  
Fefferman and Phong  established the local version for $\alpha = 1$ 
by a double commutator estimate and the theory of
pseudo\-differential operators.
The latter limits the result to operators with $C^\infty$-coefficients. 
But if the coefficients are smooth then much more is true.

\begin{cor}\label{ccom51}
If $H_0$ is a subelliptic operator of order $\gamma \in \langle0,1]$
with coefficients $c_{ij} \in C_b^\infty(\Ri^d)$ then 
$D(H^\alpha) \subseteq D(\Delta^{\alpha \gamma})$  and $(\ref{esub1.11})$  is valid     
for all $\alpha \geq 0$.
\end{cor}

Our proof of Theorem~\ref{tcom31}
 uses a double commutator estimate combined with techniques of 
functional analysis \cite{DrS} \cite{Rob7}.

Our second result deals with operators of the special form
\begin{equation}
H_0=\sum^N_{i=1}X_i^* \, X_i
\label{etcom11;2}
\end{equation}
constructed from $C_b^\infty$-vector fields $X_1,\ldots,X_N$,
i.e., vector fields on $\Ri^d$ with coefficients in $C_b^\infty(\Ri^d)$, satisfying a uniform version of 
H\"ormander's criterion for hypoellipticity.
Specifically, if $r \in \Ni$ then the vector fields $X_1,\ldots,X_N$ are defined 
to satisfy the  {\bf uniform H\"ormander condition of order $r$} if
each $C_b^\infty$-vector field $X$ can be expressed as a linear combination
\[
X=\sum_{\alpha: \; 1\leq |\alpha|\leq r}\psi_\alpha\, X_{[\alpha]}
\]
with $\psi_\alpha\in C_b^\infty(\Ri^d)$    
where $\alpha=(i_1,\ldots,i_n)$ is a multi-index with  $i_k\in\{1,\ldots,N\}$, $|\alpha|=n$,
and $X_{[\alpha]}=[X_{i_1},[X_{i_2},\ldots,[X_{i_{n-1}},X_{i_n}]\ldots]]$ is the corresponding 
multi-commutator.
This version of the H\"ormander condition was introduced by Kusuoka and Stroock
(see \cite{KuS3}, Condition (H) on page 400).
In Section~\ref{Shorf2} we present several different characterizations of the uniform H\"ormander condition.

\begin{thm} \label{thor1}
Let $H_0$ be given by $(\ref{etcom11;2})$ where
 $X_1,\ldots,X_N$ are $C_b^\infty$-vector fields on $\Ri^d$
satisfying the uniform H\"ormander condition of order $r$.
Further let $H$ denote the closure of~$H_0$. 
If $n \in \Ni$  then $D(H^n) \subseteq D(\Delta^{n/r})$ and 
there exists a $c >0$ such that
\begin{equation}
c \, \|(I+H)^n \varphi\|_2
\geq \|\Delta^{n/r} \varphi\|_2 
\label{ehor1;3}
\end{equation}
for all $\varphi\in D(H^n)$.
\end{thm}

Theorem~\ref{thor1} follows from Theorem~\ref{tcom31}, or Corollary~\ref{ccom51}, once one establishes
that the operator $H_0$ given by (\ref{etcom11;2}) satisfies the estimate (\ref{esub1.1})  with $\gamma=r^{-1}$.
The latter is a global version of H\"ormander's key estimate (\cite{Hor1}, Theorem~4.3).
H\"ormander's  argument established a local version of (\ref{esub1.1})   for all $\gamma\in\langle0,r^{-1}\rangle$.
Rothschild and Stein \cite{RS}  subsequently  established the local estimate for the optimal value  $\gamma=r^{-1}$. 
The arguments of Rothschild and Stein, which also establish optimal local versions of 
the estimates (\ref{ehor1;3}), are  based on an application of their general lifting theory.
Our arguments are completely independent of this technique and provide an alternative  proof of the optimal local results.

\section{Improvement properties} \label{Shor3a}

In this section we prove Theorem~\ref{tcom31} by use of commutator estimates.
Commutator theory was initially developed by Glimm and Jaffe \cite{GJ1} to derive self-adjointness and
regularity properties of quantum fields.
It has since developed into a useful tool for various applications in mathematical physics 
(see, for example, \cite{GJ}, Section~19.4, \cite{RS2}, Section~X.5, \cite{Far}, Section~II.12, or \cite{CFKS}, Section~4.1).
Most of these applications are based on single commutator estimates but the analysis of degenerate
operators requires double commutator  estimates \cite{DrS} \cite{Rob7}.

In the sequel we need to estimate double commutators such as 
$[\Delta,[\Delta,H_0]]$ or  analogous commutators with powers and
fractional powers of $\Delta$.
If the coefficients of $H_0$ are in $C_b^\infty(\Ri^d)$ then the commutators are defined as operators
on $W^{\infty,2}(\Ri^d)$.
If, however, the coefficients of $H_0$ are only twice differentiable  then the double commutators have to be defined as 
sesquilinear forms on  $W^{\infty,2}(\Ri^d)\times  W^{\infty,2}(\Ri^d)$.
In general, if $A,B$ are two symmetric operators in a Hilbert space $\ch$ 
and $\cd \subset D(A) \cap D(B)$ is a subspace of $\ch$ then the commutator $[B,A]$
is defined as a sesquilinear form, with form domain $\cd$, by
\[
(\psi,[B,A] \varphi)
= (B \psi, A \varphi) - (A \psi, B \varphi)
\;\;\; . 
 \]
Moreover, if $A,B_1,B_2$ are three symmetric operators in $\ch$ 
and $\cd \subset D(A) \cap D(B_2 B_1) \cap D(A B_1) \cap D(A B_2)$
then the double commutator is defined as a sesquilinear form $[B_1,[B_2,A]] $, with form domain $\cd$, by
\[
(\psi, [B_1,[B_2,A]] \varphi)
= (B_2 B_1 \psi, A \varphi) - (A B_1 \psi, B_2 \varphi) + (A \psi, B_2 B_1 \varphi) - (A B_2 \psi, B_1 \varphi)
\;\;\; .  \]
Although this is a slight abuse of notation it should not cause any confusion.
Subsequent calculations of commutators involving differential operators and multiplication operators have to be  interpreted in this form sense.
Such  commutators simplify by use of the relations 
$[\partial_i,c\,]\varphi=(\partial_i c)\varphi$ where $c$ is  a  differentiable 
function acting as a multiplication operator.

Double commutators enter estimates through the two identities
\begin{equation}
\RRe (B_2 \varphi, [B_1,A] \varphi)
= 2^{-1} (\varphi, [B_2,[B_1,A]] \varphi)
\label{elcom9;1}
\end{equation}
and 
\begin{equation}
\RRe (A\varphi,B^2\varphi)
 = (B\varphi,AB\varphi)+2^{-1}(\varphi,[B,[B,A]]\varphi)
\label{ecomm2.00}
\end{equation}
for all $\varphi \in \cd$.
In particular if $A\geq0$ the first term on the right of (\ref{ecomm2.00}) is 
positive and the double commutator gives a lower bound.

Throughout the rest of this section we set $L=I+\Delta$ and let $S_t$ denote the self-adjoint 
contraction semigroup generated by $L$.
Further we let $H_0$ be the second-order positive operator in divergence form
with coefficients $c_{ij}$ given by (\ref{esub1.0}) where the $c_{ij} \in W^{m+1,\infty}(\Ri^d)$
and $m \in \Ni$ are fixed.

\begin{lemma}\label{lcom9} 
The following  commutator estimates  are valid.
\begin{tabel}
\item\label{lcom9-1.1}
There is a $c > 0$ such that 
\[
\sum_{k=0}^d |(\psi, [\partial_k^m,[\partial_k^m,H_0]] \varphi)|
\leq c \,\|L^{m/2} \psi\|_2 \, \|L^{m/2} \varphi\|_2
\]
for all $\varphi,\psi\in W^{\infty,2}(\Ri^d)$.
\item\label{lcom9-1}
There is a $c > 0$ such that 
\[
|(\psi,[L^m,[L^m,H_0 ]]\varphi)|
\leq c \sum_{n=m}^{3m}
     \|L^{n/2}\psi\|_2 \, \|L^{(4m-n)/2}\varphi\|_2
\]
for all $\varphi,\psi\in W^{\infty,2}(\Ri^d)$.
\item\label{lcom9-2} 
If, moreover,  $c_{ij} \in W^{3,\infty}(\Ri^d)$
then there is a $c>0$ such that 
\[
|(\psi,[L,[L,H_0 ]] \varphi)|
\leq c \, \|L \psi\|_2 \, \|L \varphi\|_2
\]
for all $\varphi,\psi\in W^{\infty,2}(\Ri^d)$.
\end{tabel}
\end{lemma}
\proof\ 
The proof is by straightforward calculation using the fact that the 
coefficients are $m+1$ times differentiable.\hfill$\Box$

\ruimte

The lemma has an important corollary which is  in two parts.
The first was a  key observation of \cite{Rob7}.
The second will be used in our analysis of H\"ormander operators in Section~\ref{Shorf3}.

\begin{cor}\label{chord2.1}  
The following  commutator estimates  are valid.
\begin{tabel}
\item\label{chord2.1-1}
There is a $c > 0$ such that 
\[
|(\psi, [S_t,[S_t,H_0]]\varphi)|\leq c\,\|\psi\|_2\,\|\varphi\|_2
\]
uniformly for all  $\varphi,\psi\in W^{\infty,2}(\Ri^d)$ and $t>0$.
\item\label{chord2.1-2}
If, moreover, $c_{ij} \in W^{3,\infty}(\Ri^d)$ then there is a $c > 0$ such that 
\[
|(\psi, [S_t,[S_t,H_0]]\varphi)|\leq c\,\|(I-S_t)\psi\|_2\,\|(I-S_t)\varphi\|_2
\]
uniformly for all  $\varphi,\psi\in W^{\infty,2}(\Ri^d)$ and $t>0$.
\end{tabel}
\end{cor}
\proof\
The proof of both statements is based on the identity
\[
(\psi, [S_t, [S_t,H_0]]  \varphi) 
= \int_0^t du \int_0^t dv \, (S_{u+v} \psi, [L,[L,H_0]] S_{2t-u-v} \varphi)
\;\;\;.  \]
If $c > 0$ is as in Lemma~\ref{lcom9}.\ref{lcom9-1} applied with $m=1$, then
\begin{eqnarray*}
|(\psi, [S_t, [S_t,H_0]]  \varphi)|
& \leq & c \sum_{n=1}^3
      \int_0^t du \int_0^t dv \, 
   \|L^{n/2} S_{u+v}  \psi\|_2 \, \|L^{(4-n)/2} S_{2t-u-v} \varphi\|_2  \\
& \leq & c  \sum_{n=1}^3
 \Big( \int_0^1 du \, u^{-n/4} \, (1-u)^{-(4-n)/4} \Big)^2
  \|\psi\|_2\,  \|\varphi\|_2^2  
\end{eqnarray*}
for all $\varphi,\psi\in W^{\infty,2}(\Ri^d)$ and $t>0$ and the first statement follows.

The second statement follows by using Lemma~\ref{lcom9}.\ref{lcom9-2}
and the Cauchy--Schwarz inequality to obtain the bounds
\[
|(\psi, [S_t, [S_t,H_0]]  \varphi) |
\leq c\,\Big( \int_0^t du \int_0^t dv \,\|L S_{u+v} \psi\|_2^2 \Big)^{1/2} \cdot
\Big( \int_0^t du \int_0^t dv \,\|L S_{2t-u-v} \varphi\|_2^2 \Big)^{1/2} 
\;\;\;.
\]
Then, however, one has
\[
\int_0^t du \int_0^t dv \,\|L S_{u+v} \psi\|_2^2
\leq \int_0^t du \int_0^t dv \,\|L S_{(u+v)/2} \psi\|_2^2
= \|(I-S_t) \psi\|_2^2
\]
with a similar estimate for the $\varphi$-factor.\hfill$\Box$

\ruimte

The foregoing commutator estimates allow one to extend the  argument
used to prove Theorem~2.10 in \cite{Rob7} and to conclude the essential self-adjointness of $H_0$.
Further if  $c_{ij} \in C_b^\infty(\Ri^d)$ one can deduce from Theorems~2.16 and 2.17 of \cite{Rob7} that the  
closure $H$ of $H_0$ generates a semigroup $T$ which leaves the Sobolev spaces 
$W^{\sigma,2}(\Ri^d)=D(L^{\sigma/2})$ invariant.
We will give  shorter  self-contained proofs of these results
and establish a key  invariance property for  $c_{ij} \in W^{m+1,\infty}(\Ri^d)$.

First, note that $H_0$ maps 
$W^{\infty,2}(\Ri^d)$ into $W^{m,2}(\Ri^d)$ since  $c_{ij} \in W^{m+1,\infty}(\Ri^d)$.  
Secondly, note that for all $\sigma \geq 0$ the space $W^{\sigma,2}(\Ri^d)$ is a Hilbert space with respect to the
inner product $\langle \,\cdot\,,\,\cdot\,\rangle_\sigma$ given  by
$\langle\psi,\varphi\rangle_\sigma = (L^{\sigma/2} \psi, L^{\sigma/2} \varphi)$.
Moreover, if $n \in \Ni_0$ then the inner product $\langle \,\cdot\,,\,\cdot\,\rangle_n'$ on 
$W^{n,2}(\Ri^d)$ defined by 
\[
\langle\psi,\varphi\rangle_n'
= \sum_{k=0}^d (\partial_k^n \, \psi, \partial_k^n \, \varphi)
\]
is norm-equivalent to $\langle \,\cdot\,,\,\cdot\,\rangle_n$.

\begin{prop} \label{pcom15}
Let $H_0$ be a positive, symmetric, second-order, divergence form operator with coefficients 
$c_{ij} \in W^{m+1,\infty}(\Ri^d)$ where $m \in \Ni$.
Then the operator $H_0$ is essentially self-adjoint on $W^{\infty,2}(\Ri^d)$.
If $T$ is the self-adjoint contraction semigroup generated by the 
closure $H$ of $H_0$ then $T$ leaves the Sobolev spaces $W^{\sigma,2}(\Ri^d)$ invariant
for all $\sigma \in [0,m]$.
Moreover, the restriction of $T$ to $W^{\sigma,2}(\Ri^d)$ is a continuous  semigroup on $W^{\sigma,2}(\Ri^d)$
and  $W^{\infty,2}(\Ri^d)$ is a core for  its generator $H_{(\sigma)}$.
\end{prop}

The proof consists of verifying the criteria of the Lumer--Phillips theorem 
on  the Hilbert spaces $L_2(\Ri^d)$ and 
$W^{m,2}(\Ri^d)$ with inner product $\langle \,\cdot\,,\,\cdot\,\rangle_m'$.

\begin{lemma}\label{lcomd2.1}
There is an $\omega\geq0$ such that 
\begin{equation}
\RRe\langle \varphi, (H_0+\omega I)) \varphi \rangle_m' \geq 0
\label{ecomd2.1}
\end{equation}
for all $\varphi\in W^{\infty,2}(\Ri^d)$.
\end{lemma}
\proof\
Since $(i \partial_k)^m$ is symmetric one deduces from (\ref{ecomm2.00}) that
\begin{eqnarray*}
\RRe \langle\varphi, H_0 \varphi \rangle_m'
& = & \sum_{k=0}^d (\partial_k^m \varphi, \partial_k^m H_0 \varphi)  \\
& = & \sum_{k=0}^d (\partial_k^m \varphi, H_0 \partial_k^m \varphi) 
        + (-1)^m 2^{-1} (\varphi, [\partial_k^m, [\partial_k^m, H_0]] \varphi)  \\
& \geq & - 2^{-1} c \,\|L^{m/2} \varphi\|_2^2
= - 2^{-1} c \, \langle\varphi, \varphi \rangle_m
\end{eqnarray*}
for all $\varphi \in W^{\infty,2}(\Ri^d)$,
where $c$ is the constant in Lemma~\ref{lcom9}.\ref{lcom9-1.1}.
Then (\ref{ecomd2.1}) is valid because the norms associated with the inner products
$\langle \,\cdot\,,\,\cdot\,\rangle_m$ and $\langle \,\cdot\,,\,\cdot\,\rangle_m'$
are equivalent.\hfill$\Box$

\begin{lemma}\label{lcomd2.2}
There exists an $\varepsilon > 0$ such that 
$(I + \varepsilon H_0) W^{\infty,2}(\Ri^d)$ is dense in $L_2(\Ri^d)$ and $W^{m,2}(\Ri^d)$.
\end{lemma}
\proof\
Let $n \in \{ 0,m \} $.
We establish  below that there exists a $c > 0$ such that 
\begin{equation}
- \RRe \langle \varphi, H_0 S_{2t} \varphi \rangle_n'
\leq c\, \|L^{n/2} \varphi\|_2^2
\label{epcom15;1}
\end{equation}
uniformly for all $t > 0$ and $\varphi \in W^{\infty,2}(\Ri^d)$.
It then follows by continuity that (\ref{epcom15;1}) is valid 
uniformly for all $t > 0$ and $\varphi \in D(L^{n/2})$.
Moreover, there exists a $c_1 > 0$ such that 
\[
\|L^{n/2} \varphi\|_2^2
\leq c_1 \sum_{k=0}^d \|\partial_k^n \varphi\|_2^2
\]
for all $\varphi \in W^{\infty,2}(\Ri^d)$.
Now set $\varepsilon = (2c c_1)^{-1}$.
Let $\varphi \in D(L^{n/2})$ and suppose that the inner product
$\langle \varphi, (I + \varepsilon H_0) \psi \rangle_n' = 0$ 
for all $\psi \in W^{\infty,2}(\Ri^d)$.
Then $S_{2t} \varphi \in W^{\infty,2}(\Ri^d)$ and 
\begin{eqnarray*}
c_1^{-1} \|S_t L^{n/2} \varphi\|_2^2
\leq \sum_{k=0}^d \|S_t \partial_k^n \varphi\|_2^2
& = & \langle \varphi, S_{2t} \varphi \rangle_n'  \\[0pt]
& = & - \varepsilon \RRe \langle \varphi, H_0 S_{2t} \varphi \rangle_n'  
\leq c \, \varepsilon \|L^{n/2} \varphi\|_2^2
= (2 c_1)^{-1} \|L^{n/2} \varphi\|_2^2
\end{eqnarray*}
for all $t > 0$.
So 
$
\|L^{n/2} \varphi\|_2^2
= \lim_{t \downarrow 0} \|S_t L^{n/2} \varphi\|_2^2
\leq 2^{-1} \|L^{n/2} \varphi\|_2^2
$
and $\varphi = 0$.
Therefore it remains to prove (\ref{epcom15;1}).

Let $t > 0$ and $\varphi \in W^{\infty,2}(\Ri^d)$.
The starting point is the identity
\begin{equation}
\RRe \langle \varphi, H_0 S_{2t} \varphi \rangle_n'
= \sum_{k=0}^d \RRe (\partial_k^n \varphi, H_0 S_{2t} \partial_k^n \varphi)
   + \sum_{k=0}^d \RRe (\partial_k^n \varphi,[\partial_k^n, H_0 ] S_{2t}\varphi)  
\;\;\; .
\label{epcom15;2}
\end{equation}
We will bound the two terms separately.

The first term satisfies the identity
\[
\sum_{k=0}^d \RRe (\partial_k^n \varphi, H_0 S_{2t} \partial_k^n \varphi)
= \sum_{k=0}^d ( S_t \partial_k^n \varphi, H_0  S_t \partial_k^n \varphi) 
   + \sum_{k=0}^d 2^{-1} (\partial_k^n \varphi,[S_t,[S_t,H_0]] \partial_k^n \varphi)
\]
where we have again used (\ref{ecomm2.00}).
Therefore if $c > 0$ is as in Corollary~\ref{chord2.1}.\ref{chord2.1-1} then 
\[
\sum_{k=0}^d \RRe \, (\partial_k^n \varphi, H_0 S_{2t} \partial_k^n \varphi)
\geq -2^{-1} c \,  \|L^{n/2} \varphi\|_2^2
\;\;\; .  \]
Note that $c$ is independent of $t$ and $\varphi$.
If $n=0$ this completes the proof since the second term in (\ref{epcom15;2}) is identically zero.

If $n = m$ then, by (\ref{elcom9;1}), the second term in (\ref{epcom15;2}) satisfies the identity
\begin{eqnarray}
\sum_{k=0}^d \RRe \,( \partial_k^m \varphi,[\partial_k^m, H_0 ] S_{2t}\varphi)  
& = & (-1)^m \sum_{k=0}^d 2^{-1}(S_t\varphi,[\partial_k^m,[\partial_k^m,H_0]]S_t\varphi)  \nonumber  \\*[-5pt]
& & \hspace{20mm} {}
   - \sum_{k=0}^d \RRe \, ( \partial_k^m \varphi,[S_t,[\partial_k^m, H_0 ]] S_t\varphi)  
\;\;\; . 
\label{epcom15;2.1}
\end{eqnarray}
But the first term in (\ref{epcom15;2.1}) satisfies
\begin{eqnarray*}
\sum_{k=0}^d 2^{-1}  |(S_t \varphi, [\partial_k^m, [\partial_k^m,H_0]] S_t \varphi) | 
\leq  2^{-1} c \, \|L^{m/2} S_t \varphi\|_2^2
\leq 2^{-1} c\, \|L^{m/2} \varphi\|_2^2
\end{eqnarray*}
with $c$ the constant in Lemma~\ref{lcom9}.\ref{lcom9-1.1}.
Finally, if $k \in \{ 1,\ldots,d \} $ then we estimate the last term as follows.
The Cauchy inequality gives
\[
|(\partial_k^m \varphi, [S_t , [\partial_k^m , H_0]] S_t \varphi)|
\leq \|L^{m/2} \varphi\|_2 \, \|[S_t , [\partial_k^m , H_0]] S_t \varphi\|_2
\;\;\; .  \]
Moreover, 
\[
[\partial_k^m , H_0]
= \sum_{i,j=0}^d \sum_{p=0}^{m-1} {m \choose p+1} 
           \Big( (\partial_i \partial_k^{p+1} c_{ij}) \partial_k^{m-p-1} \partial_j 
    + (\partial_k^{p+1} c_{ij}) \partial_i \partial_k^{m-p-1} \partial_j \Big)
\;\;\; .  \]
Therefore 
\begin{eqnarray*}
\|[S_t , [\partial_k^m , H_0]] S_t \varphi\|_2
& \leq & \sum_{i,j=1}^d \sum_{p=0}^{m-1} {m \choose p+1} 
    \Big( 
   \|[S_t, (\partial_i \partial_k^{p+1} c_{ij})] S_t\|_{2 \to 2}   \\[-5pt]
& & \hspace{50mm} {}
    + \|[S_t, (\partial_k^{p+1} c_{ij})]\,\partial_i S_t\|_{2 \to 2} \Big)  \, \|L^{m/2} \varphi\|_2
    \;\;\;.
\end{eqnarray*}
The first term on the right is clearly bounded by a multiple of $\|L^{m/2} \varphi\|_2$.
But the second satisfies a similar bound since
\begin{eqnarray*}
\|[S_t, (\partial_k^{p+1} c_{ij})]\,\partial_i S_t\|_{2 \to 2} 
&\leq & \sum^d_{l=1}\int^t_0du \,
\|\partial_lS_u (\partial_l \partial_k^{p+1} c_{ij}) \partial_i S_{2t-u}
   + S_u (\partial_l \partial_k^{p+1} c_{ij}) \partial_l \partial_i S_{2t-u}\|_{2 \to 2} \\
&\leq& \sum_{l=1}^d \|\partial_l \partial_k^{p+1} c_{ij}\|_\infty 
      \int^1_0du \Big( u^{-1/2}(2-u)^{-1/2}+(2-u)^{-1} \Big)
\end{eqnarray*}
for all $p \in \{ 0,\ldots,m-1 \} $.
This completes the proof of (\ref{epcom15;1}) and the proof of the lemma.\hfill$\Box$

\ruimte

\noindent{\bf Proof of Proposition~\ref{pcom15}\hspace{5pt}\ } 
The operator $H_0$ is positive and symmetric
on $L_2(\Ri^d)$.
It then follows from Lemma~\ref{lcomd2.2}  that it is essentially self-adjoint in $L_2(\Ri^d)$.
Therefore the self-adjoint closure $H$  generates a self--adjoint contraction semigroup
$T$  on $L_2(\Ri^d)$.
It follows from  Lemmas~\ref{lcomd2.1}, \ref{lcomd2.2}
and the Lumer--Phillips theorem, \cite{LuP} Theorem~3.1, that the operator $H_0$ is closable on $W^{m,2}(\Ri^d)$
and that its closure generates a continuous quasi-contraction semigroup on $W^{m,2}(\Ri^d)$
if $W^{m,2}(\Ri^d)$ is equipped with  the norm induced from the inner product $\langle \,\cdot\,,\,\cdot\,\rangle_m'$.
But this quasi-contraction semigroup is automatically the restriction of $T$ to $W^{m,2}(\Ri^d)$.
Moreover, it is a continuous semigroup on $W^{m,2}(\Ri^d)$ equipped with the norm induced by 
$\langle \,\cdot\,,\,\cdot\,\rangle_m$.
If $\sigma \in \langle0,m \rangle$ it follows by interpolation that $T$ leaves the Sobolev space 
$W^{\sigma,2}(\Ri^d)$ invariant and the restriction of $T$ to $W^{\sigma,2}(\Ri^d)$ is a 
continuous  semigroup on $W^{\sigma,2}(\Ri^d)$.
Let $H_{(\sigma)}$ denotes the generator on $W^{\sigma,2}(\Ri^d)$.
Then $(\lambda I + H_0) W^{\infty,2}(\Ri^d)$ is dense in $W^{m,2}(\Ri^d)$ for large $\lambda > 0$ 
by Lemma~\ref{lcomd2.2}, so $(\lambda I + H_{(\sigma)}) W^{\infty,2}(\Ri^d) = (\lambda I + H_0) W^{\infty,2}(\Ri^d)$
is dense in $W^{\sigma,2}(\Ri^d)$.
Therefore $W^{\infty,2}(\Ri^d)$ is a core for $H_{(\sigma)}$.\hfill$\Box$

\ruimte

Next we turn to the problem of  improving  order properties.
Note that if $A$ and $B$ are self-adjoint operators and $A\geq B^2$ then it is not true in general that $A^2\geq B^4$
although it is true if $A$ and $B$ commute.
The next lemma draws a similar conclusion from a double commutator bound.

\begin{lemma} \label{lcom2}
Let $\cd$ be a subspace of a Hilbert space $\ch$ and $A,B$ a symmetric and self-adjoint operator
on $\ch$,
respectively,
such that $\cd \subset D(A) \cap D(B)$ and $B \cd \subset \cd$.
Assume 
\begin{equation}
(\varphi, A \varphi) \geq \|B \varphi\|^2
\label{elcom2;1}
\end{equation}
for all $\varphi \in \cd$.
If there are  $\varepsilon \in [0,1\rangle$ and $c > 0$ such that 
\begin{equation}
|(\varphi,[B,[B,A]]\varphi)|
\leq  \varepsilon \, \|B^2\varphi\|^2 + c \, \|\varphi\|^2
\label{elcom2;2}
\end{equation}
for all $\varphi \in \cd$, then 
\[
\|A \varphi\|^2 \geq (1-\varepsilon) \|B^2 \varphi\|^2 - c \, \|\varphi\|
\]
and in particular,
\[
\|A \varphi\| \geq (1-\varepsilon) \|B^2 \varphi\| - c^{1/2} \, \|\varphi\|
\]
for all $\varphi \in \cd$.
\end{lemma} 
\proof\
One estimates  that 
\begin{eqnarray*}
\|A\varphi\|^2 + \|B^2\varphi\|^2
& \geq & 2 \RRe (A \varphi, B^2 \varphi)   \\
& = & 2 \, (B\varphi,AB\varphi)+(\varphi,[B,[B,A]]\varphi)\\
&\geq& (2 - \varepsilon) \, \|B^2\varphi\|^2 -c \, \|\varphi\|^2
\end{eqnarray*}
for all $\varphi\in\cd$ where we have successively used (\ref{ecomm2.00}), (\ref{elcom2;1})
and (\ref{elcom2;2}).
The statement of the lemma follows immediately.\hfill$\Box$

\ruimte

The double commutator estimate (\ref{elcom2;2}) is a rather weak requirement
for second-order differential operators.
For example, if $B = L$ and $A = H_0$ then Lemma~\ref{lcom9}.\ref{lcom9-2} gives the much stronger bound
\[
|(\varphi,[B,[B,A]]\varphi)|
\leq c \, \|B\varphi\|^2 
\;\;\; .  \]
But our proof of the improvement of subelliptic properties follows from application
of Lemma~\ref{lcom2} with $B$ a fractional power of $L$ and this leads to a slight 
`loss of derivatives'.
Recall that we assume $c_{ij} \in W^{m+1,\infty}(\Ri^d)$ with $m \in \Ni$.
 
\begin{lemma}\label{lcom1.1f}
For all $\rho \in [0,m\rangle$ and $\delta > 0 \vee (\rho - 2^{-1} m)$ there is a $c > 0$ such that
\[
|(\varphi,[L^\rho,[L^\rho,H_0 ]]\varphi)|
\leq c\,\|L^{\rho + \delta}\varphi\|_2^2
\]
for all $\varphi\in W^{\infty,2}(\Ri^d)$.
\end{lemma}
\proof\
The case $\rho = 0$ is trivial, so we may assume that $\rho > 0$.
Set $\tau = m^{-1} \rho \in \langle0,1\rangle$.
Then
\[
L^\rho 
= (L^m)^\tau
= c_1 \int^\infty_0d\lambda\,\lambda^{-1+\tau}L(\lambda I+L)^{-1}
\]
where $c_1 = \int^\infty_0d\lambda\,\lambda^{-1+\tau}(1+\lambda)^{-1}$.
Let $c > 0$ be as in Lemma~\ref{lcom9}.\ref{lcom9-1}.

Let $\varphi \in W^{\infty,2}(\Ri^d)$.
Then 
\begin{eqnarray*}
(\varphi, [L^\rho,[L^\rho,H_0]]\varphi)
& = & c_1^{-2} \int^\infty_0 d\lambda \int^\infty_0 d\mu \,
    (\lambda\mu)^\tau (R_\lambda R_\mu\varphi, [L^m,[L^m,H_0]] R_\lambda R_\mu\varphi)  \\
& \leq & c \, c_1^{-2} \sum_{n=m}^{3m} \int^\infty_0 d\lambda \int^\infty_0 d\mu \,
    (\lambda\mu)^\tau \|L^{n/2} R_\lambda R_\mu \varphi\|_2 \, \|L^{(4m-n)/2} R_\lambda R_\mu \varphi\|_2
\end{eqnarray*}
where $R_\lambda = (\lambda I + L^m)^{-1}$.
It follows from spectral theory that 
$\|L^\alpha R_\lambda\|_{2 \to 2} \leq (1 + \lambda)^{-(m-\alpha)/m}$ for all $\lambda > 0$ and 
$\alpha \in [0,m]$.

Let $n \in \{ m,\ldots,3m \} $. 
Set $\eta_1 = (2^{-1} n) \wedge (\rho + \delta)$ and 
$\eta_2 = (2^{-1} (4m-n)) \wedge (\rho + \delta)$.
Then 
\[
\!\|L^{n/2} R_\lambda R_\mu \varphi\|_2
\leq \|L^{(n - 2 \eta_1)/2} R_\lambda R_\mu\|_{2 \to 2}  \|L^{\eta_1} \varphi\|_2
\leq \Big( (1+\lambda) (1 + \mu) \Big)^{- (4m)^{-1} (4m-n+2\eta_1) } \|L^{\rho + \delta} \varphi\|_2
\]
for all $\lambda,\mu > 0$.
Similarly, 
\[
\|L^{(4m-n)/2} R_\lambda R_\mu \varphi\|_2
\leq \Big( (1+\lambda) (1 + \mu) \Big)^{- (4m)^{-1} (n+2\eta_2) } \|L^{\rho + \delta} \varphi\|_2
\;\;\; .  \]
Therefore 
\begin{eqnarray*}
\lefteqn{
\int^\infty_0 d\lambda \int^\infty_0 d\mu \,
    (\lambda\mu)^\tau \|L^{n/2} R_\lambda R_\mu \varphi\|_2 \, \|L^{(4m-n)/2} R_\lambda R_\mu \varphi\|_2
} \hspace{60mm} \\
& \leq & \Big( \int^\infty_0 d\lambda \, \lambda^\tau \, (1+\lambda)^{- (2m)^{-1} (2m + \eta_1 + \eta_2)} \Big)^2
     \|L^{\rho + \delta} \varphi\|_2^2
\;\;\; .
\end{eqnarray*}
So it remains to verify that the integral is finite, i.e., we have to show that 
$\eta_1 + \eta_2 > 2 m \tau = 2 \rho$.

If $\rho + \delta \leq 2^{-1} n \wedge 2^{-1} (4m-n)$ then
$\eta_1 + \eta_2 = 2 (\rho + \delta) > 2 \rho$.
If $2^{-1} n \vee 2^{-1} (4m-n) \leq \rho + \delta$ then 
$\eta_1 + \eta_2 = 2m > 2 \rho$.
Finally, if $2^{-1} n \leq \rho + \delta \leq  2^{-1} (4m-n)$ then 
$\eta_1 + \eta_2 = 2^{-1} n + \rho + \delta \geq 2^{-1} m + \rho + \delta > 2 \rho$
since $\delta > \rho - 2^{-1} m$.
Similarly $\eta_1 + \eta_2 > 2 \rho$ if $2^{-1} (4m-n) \leq \rho + \delta \leq  2^{-1} n$.
This proves the lemma.\hfill$\Box$

\ruimte

The previous lemmas can be  applied to establish an improvement of Theorem~\ref{tcom31}.

\begin{thm} \label{tcom32}
Let $H_0$ be a subelliptic operator of order $\gamma \in \langle0,1]$
with coefficients $c_{ij} \in W^{m+1,\infty}(\Ri^d)$, where  $m \in \Ni$,
and with self-adjoint closure $H$.
Further, let $\sigma \in [0,2^{-1} m \rangle$.
If $\varphi \in D(H)$ and $H \varphi \in W^{2\sigma,2}(\Ri^d)$ then 
$\varphi \in W^{2\sigma + 2 \gamma,2}(\Ri^d)$.
Moreover, there exist  $c,\omega_0 > 0$ such that 
\begin{equation}
c \, \|L^\sigma (\omega_0 I + H_{(2 \sigma)}) \varphi\|_2
\geq \|L^{\sigma + \gamma} \varphi\|_2
\label{etcom32;1}
\end{equation}
for all $\varphi \in D(H_{(2 \sigma)})$, where $H_{(2 \sigma)}$ denotes the closure of $H_0$ on $W^{2\sigma,2}(\Ri^d)$.
\end{thm}
\proof\
Since the restriction of $S$ to $W^{2\sigma,2}(\Ri^d)$ is a continuous semigroup there exists an $\omega_0 > 1$
such that $\|L^\sigma \varphi\|_2 \leq \|L^\sigma (\omega_0 I + H_0) \varphi\|_2$
for all $\varphi \in W^{\infty,2}(\Ri^d)$.
Set $\tau = 2^{-1} \gamma$.
Since $\|L^\tau \varphi\|_2^2 \leq 4 \,\|\Delta^\tau\varphi\|_2^2+\|\varphi\|_2^2$
there exists by subellipticity a $c > 0$ such that 
\[
c\,(\varphi, (\omega_0 I + H_0) \varphi)\geq \|L^\tau\varphi\|_2^2
\]
for all $\varphi \in W^{\infty,2}(\Ri^d)$.
Set $A = c \, L^\sigma (\omega_0 I + H_0) L^\sigma$ and $B = L^{\sigma + \tau}$.
Then $(\varphi, A \varphi) \geq \|B \varphi\|_2^2$ for all $\varphi \in W^{\infty,2}(\Ri^d)$.
So (\ref{elcom2;1}) in Lemma~\ref{lcom2} is satisfied with $\cd=W^{\infty,2}(\Ri^d)$.
Fix $\delta \in \langle (\sigma - 2^{-1} m + \tau) \vee 0, \tau \rangle$.
Then it follows from Lemma~\ref{lcom1.1f} that there are $c',\omega > 0$ such that 
\begin{eqnarray*}
|(\varphi,[B,[B,A]]\varphi)|
& = & c \, |(L^\sigma \varphi, [L^{\sigma + \tau}, [L^{\sigma + \tau}, H_0]] L^\sigma \varphi)|  \\
& \leq & c' \|L^{2 \sigma + \tau + \delta} \varphi\|_2^2
\leq 2^{-1} \|L^{2 \sigma + 2\tau} \varphi\|_2^2 + \omega \, \|\varphi\|_2^2
\end{eqnarray*}
for all $\varphi \in W^{\infty,2}(\Ri^d)$.

So by Lemma~\ref{lcom2}
\[
c^2 \, \|L^\sigma (\omega_0 I + H_0) L^\sigma \varphi\|_2^2
\geq 2^{-1} \|L^{2 \sigma + 2\tau} \varphi\|_2^2 - \omega \,\|\varphi\|_2^2
\geq 2^{-1} \|L^{2 \sigma + 2\tau} \varphi\|_2^2 - \omega\, \|L^{2 \sigma} \varphi\|_2^2
\]
for all $\varphi \in W^{\infty,2}(\Ri^d)$.
Since $L^\sigma$ is a bijection from $W^{\infty,2}(\Ri^d)$ onto $W^{\infty,2}(\Ri^d)$ one may 
replace $L^\sigma\varphi$ by $\varphi$  and then (\ref{etcom32;1})
follows by rearrangement uniformly for all $\varphi \in W^{\infty,2}(\Ri^d)$.
But $W^{\infty,2}(\Ri^d)$ is a core for $H_{(2 \sigma)}$ by Proposition~\ref{pcom15}.
So (\ref{etcom32;1}) is valid for all $\varphi \in D(H_{(2 \sigma)})$.

The first statement follows immediately from the second.\hfill$\Box$

\ruimte

\noindent{\bf Proof of Theorem~\ref{tcom31}\hspace{5pt}\ }
This follows immediately from Theorem~\ref{tcom32} by interpolation and a telescopic argument.\hfill$\Box$

\ruimte

We conclude this section with four remarks on Theorem~\ref{tcom31}.

First, 
Theorem~\ref{tcom31} can be rephrased in terms of order relations.
The estimate (\ref{etcom32;1}) with $\sigma=0$ is equivalent to the quadratic form estimate
\[
c^2\,(I+H)^2\geq L^{2\gamma}
\;\;\;.
\]
Then since the order relation between positive self-adjoint operators is respected by taking fractional 
powers one has
\[
c^{2\alpha}\,(I+H)^{2\alpha}\geq L^{2\alpha \gamma}
\]
for all $\alpha\in\langle0,1]$.
But (\ref{etcom32;1}) with $\sigma>0$ is equivalent to the estimate
\[
c^2\,(I+H)L^{2\sigma} (I+H)\geq L^{2\sigma+2\gamma}
\]
and then the previous argument can be iterated to obtain the order relations covered by 
Theorem~\ref{tcom31}.

Secondly, 
if $c_{ij}\in W^{2,\infty}(\Ri^d)$ then Theorem~\ref{tcom31} establishes that the subellipticity 
condition (\ref{esub1.1}) implies the estimate (\ref{esub1.11}) with $\alpha=1$.
But the foregoing observation on order properties establishes the converse.
Thus (\ref{esub1.1}) is equivalent to  (\ref{esub1.11}) with $\alpha=1$.
This is a global strengthening of the first statement in Theorem 1 of Fefferman and Phong \cite{FP}.

Thirdly, 
the statement of Theorem~\ref{tcom31}  is partly redundant  since
 the closed graph theorem implies that if  the inclusion 
$D(H^\alpha) \subseteq D(\Delta^{\alpha \gamma})$  is valid for one $\alpha\geq0$ then 
there exists a $c > 0$ such that (\ref{esub1.11}) is valid.

Finally, 
if $\gamma=1$, i.e., if $H_0$ is strongly elliptic, then the statement of the theorem is also valid 
for  $\alpha = 2^{-1}(m + 1 + \gamma^{-1}) = 2^{-1}(m+2)$,
and this is  the best estimate one could expect for operators with coefficients $c_{ij}\in W^{m+1,\infty}(\Ri^d)$.
The extension is a consequence of the theorem, applied with $\alpha = 2^{-1}(m+1)$, together 
with a simple commutator estimate.
In fact  for $\gamma = 1$ the domain inclusion in the theorem is an equality
which is also valid on the $L_p$-spaces if $p \in \langle1,\infty\rangle$ 
(see \cite{ER19}, Theorem 1.5.II).
 
\section{$C_b^\infty$-flows}\label{Shor3}

In this section we prepare  the discussion of elliptic operators (\ref{etcom11;2}) of H\"ormander type 
by recalling  some basic properties of the flows corresponding 
to $C_b^\infty$-vector fields. 
We also  give several estimates for products and 
commutators of such flows.
Local estimates of a similar nature are an important feature in the work of 
H\"ormander \cite{Hor1} and Nagel, Stein and Wainger \cite{NSW} but our emphasis is
on  estimates which are   uniform over $\Ri^d$.
The uniform H\"ormander condition is not relevant in this section.

Let $X$ be a $C_b^\infty$-vector field on $\Ri^d$ with coefficients $a_i$.
Then it follows from the theory of ordinary differential equations 
that there exists a unique $C^\infty$-function $f \colon \Ri \times \Ri^d \to \Ri^d$
such that 
\[
\frac{\partial f_i}{\partial t} (t,x) = a_i(f(t,x))
\;\;\;\; \mbox{and} \;\;\;\; 
f(0,x) = x
\]
for all $x \in \Ri^d$, $t \in \Ri$ and $i \in \{ 1,\ldots,d \} $.
We adopt the conventional notation $\exp(tX)(x) = f(t,x)$.
Then for all $\varphi \in C^\infty(\Ri^d)$ and $t \in \Ri$ we define 
$e^{tX} \varphi \in C^\infty(\Ri^d)$ by 
$(e^{tX} \varphi)(x) = \varphi(\exp(tX)(x))$.
The relevant properties of these maps  are summarized as follows.

\begin{lemma} \label{lhor3}
Let $X$ be a $C_b^\infty$-vector field on $\Ri^d$.
Then one has the following.
\begin{tabel}
\item \label{lhor3-1}
$\exp(tX)(\exp(sX)(x)) = \exp((t+s)X)(x)$ for all $x \in \Ri^d$ and $t,s \in \Ri$.
Hence for each $t \in \Ri$ the map $\exp(tX)$ is a diffeomorphism of $\Ri^d$.
\item \label{lhor3-2}
If $\varphi \in C^\infty(\Ri^d)$ and $t \in \Ri$ then 
\[
\varphi(\exp(tX)(x)) 
= ( e^{tX} \varphi )(x)
\sim \sum_{n=0}^\infty t^n \, n!^{-1} (X^n \varphi)(x)
\]
for all $x \in \Ri^d$, where $\sim$ denotes the Taylor series $($in $t)$ around $0$.
\end{tabel}
\end{lemma}

We also need some quantitative estimates.
It is convenient to introduce a multi-index notation.
For all $N \in \Ni$ and $n \in \Ni_0$ set
\[ 
J_n(N) = \bigoplus_{k=0}^n \{ 1,\ldots,N \} ^k 
\;\;\;\;\;,
\;\;\;\;\;
J(N) = \bigoplus_{k=0}^\infty \{ 1,\ldots,N \} ^k
\]
and let $J^+_n(N)$, $J^+(N)$ denote the corresponding sets 
with the restrictions $k \geq 1$.

One can prove the next lemma with the aid of Gronwall's lemma and induction.

\begin{lemma} \label{lhor3.1}
Let $X$ be a $C_b^\infty$-vector field.
\begin{tabel}
\item \label{lhor3.1-2}
For all $k \in \Ni$ there exists an $M > 0$ such that 
\[
|\partial_t^k \, \exp(tX)(x)| 
\leq M 
\]
uniformly for all $t \in \Ri$ and $x \in \Ri^d$.
\item \label{lhor3.1-3}
For all $\alpha \in J(d)$ and $k \in \Ni_0$ with $|\alpha| + k \geq 1$
there exist
$M,\omega > 0$ such that 
\[
|\partial_t^k \, \partial_x^\alpha \exp(tX)(x)| 
\leq M e^{\omega |t|}
\]
uniformly for all $t \in \Ri$ and $x \in \Ri^d$.
\item \label{lhor3.1-4}
There are $M,\omega > 0$ such that 
\[
\|e^{tX} \varphi\|_2 \leq M e^{\omega |t|} \, \|\varphi\|_2
\]
for all $t \in \Ri$ and $\varphi \in C_c^\infty(\Ri^d)$.
\end{tabel}
\end{lemma}

Next we need several estimates which follow from the 
Campbell--Baker--Hausdorff formula.
The first of these is an estimate for the product of two flows generated by
$C_b^\infty$-vector fields.
The key observation is contained in the following lemma. 

\begin{lemma} \label{lhor6}
Let $Y_1$ and $Y_2$ be $C_b^\infty$-vector fields
and let $N \in \Ni \backslash \{ 1 \} $.
Then there exist $Z_2,\ldots,Z_N$ with 
$Z_j \in \spann \{ Y_{[\alpha]} : \alpha \in J(2), \; |\alpha| = j \}$
for all $j \in \{ 2,\ldots,N \} $, such that 
for all $\varphi \in C^\infty(\Ri^d)$ and $x \in \Ri^d$ one has
\[
\bigg( \frac{d}{dt} \bigg)^k
   \varphi( \exp(t(Y_1+Y_2)) \exp(-t Y_1) \exp(-t Y_2) \exp(-t^2 Z_2) \ldots
            \exp(- t^N Z_N)(x) ) \bigg|_{t=0}
= 0
\]
for all $k \in \{ 1,\ldots,N \} $.
\end{lemma}
\proof\
This follows from the Campbell--Baker--Hausdorff formula as in the discussion
preceding Lemma~4.5 of \cite{Hor1}. 
See in particular \cite{Hor1}, pp.\ 160--161.\hfill$\Box$

\ruimte
As a direct consequence one has the following estimate which is uniform over $\Ri^d$.

\begin{prop} \label{phor301}
Let $Y_1$ and $Y_2$ be $C_b^\infty$-vector fields
and let $N \in \Ni \backslash \{ 1 \} $.
Then there exist $c>0$ and $Z_2,\ldots,Z_N$ with 
$Z_j \in \spann \{ Y_{[\alpha]} : \alpha \in J(2), \; |\alpha| = j \}$
for all $j \in \{ 2,\ldots,N \} $, such that 
\[
|\exp( t(Y_1 + Y_2) ) \exp(-t Y_1) \exp(-t Y_2) \exp(-t^2 Z_2) \ldots
            \exp(- t^N Z_N)(x) -x| 
 \leq c\,t^{N+1}
\]
uniformly for all $x\in\Ri^d$ and $t\in[-1,1]$.
\end{prop}
\proof\
Define $\Phi \colon \Ri^d \times \Ri \to \Ri^d$ by
\[
\Phi(x,t) 
= \exp( t(Y_1 + Y_2) ) \exp(-t Y_1) \exp(-t Y_2) \exp(-t^2 Z_2) \ldots
            \exp(- t^N Z_N)(x) 
\;\;\; .  \]
If $\varphi \in C^\infty(\Ri^d)$ then it follows from Lemma~\ref{lhor6} and
the Taylor integral  remainder formula that 
\begin{eqnarray*}
|\varphi(\Phi(x,t)) - \varphi(x)|
& = & |\varphi(\Phi(x,t)) - \varphi(\Phi(x,0))|  \\
& = & \Big|  N!^{-1} \int_0^t ds \,(t-s)^N \, \partial_s^{N+1} \varphi(\Phi(x,s)) \Big|
\end{eqnarray*}
for all $x \in \Ri^d$ and $t \in \Ri$.
Now apply the above to $\varphi = \pi_k$.
It follows from Lemma~\ref{lhor3.1}.\ref{lhor3.1-3} that there exists an $M > 0$ 
such that 
$|\partial_s^{N+1} \pi_k(\Phi(x,s))| \leq M$ uniformly for all $x \in \Ri^d$ and $s \in [-1,1]$.
Then
\[
|\pi_k(\Phi(x,t)) - x_k|
\leq  N!^{-1} \int_0^t ds \,  (t-s)^N \, M 
= M \, (N+1)!^{-1} \, t^{(N+1)}
\]
for all $x \in \Ri^d$, $t \in [-1,1]$ and $k \in \{ 1,\ldots,d \} $.
So $|\Phi(x,t) - x| \leq M d \, t^{N+1}$
for all $x \in \Ri^d$ and $t \in [-1,1]$.\hfill$\Box$

\ruimte

Finally we give an estimate comparing the flow
 generated by a combination of 
multi-commutators in terms of products of the elementary flows.

\begin{lemma} \label{lhor848}
Let $X_1,\ldots,X_N$ be $C_b^\infty$-vector fields and $s \in \Ni$.
Then there exists an $M > 0$ such that for all $b \colon J_s^+(N) \to [-1,1]$
there are $n \in \{ 1,\ldots,3(2d)^s\} $, $i_1,\ldots,i_n \in \{ 1,\ldots,N \} $
and $a_1,\ldots,a_n \in [-M,M]$ such that 
\[
\bigg( \frac{d}{dt} \bigg)^k
   \varphi( \exp( \sum_{\alpha \in J_s^+(N)} b(\alpha) \, t^{|\alpha|} X_{[\alpha]})
        \exp(-a_1 t X_{i_1}) \ldots \exp(-a_n t X_{i_n})(x) ) \bigg|_{t=0}
= 0
\]
for all $k \in \{ 1,\ldots,N \} $, $\varphi \in C^\infty(\Ri^d)$ and $x \in \Ri^d$.
\end{lemma}
\proof\
This follows from the arguments in the proof of Lemma~2.22 in \cite{NSW}.\hfill$\Box$

\begin{prop} \label{phor843}
Let $s \in \Ni$ and $X_1,\ldots,X_N$ be $C_b^\infty$-vector fields.
Then there exist  $M,M' > 0$ such that for all $\delta \in \langle0,1]$
and $b \colon J_s^+(N) \to [-1,1]$ with $|b(\alpha)| \leq \delta^{|\alpha|}$
for all $\alpha \in J_s^+(N)$ there are $n \in \{ 1,\ldots,3(2d)^s \} $,
$i_1,\ldots,i_n \in \{ 1,\ldots,N \} $
and $a_1,\ldots,a_n \in [-M \delta,M \delta]$ such that 
\[
|\exp( \sum_{\alpha \in J_s^+(N)} b(\alpha) \, t^{\alpha} \, X_{[\alpha]})(x)
    \exp(-a_1 t X_{i_1}) \ldots \exp(-a_n t X_{i_n})(x) - x| 
\leq M' \, (\delta \, |t|)^{s+1}
\]
uniformly for all $x \in \Ri^d$ and $t \in [-1,1]$.
\end{prop}
\proof\
Let $M > 0$ be as in Lemma~\ref{lhor848}.
Let $\delta \in \langle0,1]$ and $b \colon J_s^+(N) \to [-1,1]$ with 
$|b(\alpha)| \leq \delta^{|\alpha|}$
for all $\alpha \in J_s^+(N)$.
Then $\delta^{-|\alpha|} \,| b(\alpha)| \leq 1$ for all $\alpha \in J_s^+(N)$.
By Lemma~\ref{lhor848}
there are $n \in \{ 1,\ldots,3(2d)^s \} $, $i_1,\ldots,i_n \in \{ 1,\ldots,N \} $
and $a_1,\ldots,a_n \in [-M,M]$ such that 
\[
\bigg( \frac{d}{dt} \bigg)^k
   \varphi( \exp( \sum_{\alpha \in J_s^+(N)} \delta^{-|\alpha|} b(\alpha) \, t^{|\alpha|} X_{[\alpha]})
        \exp(-a_1 t X_{i_1}) \ldots \exp(-a_n t X_{i_n})(x) ) \bigg|_{t=0}
= 0
\]
for all $k \in \{ 1,\ldots,N \} $, $\varphi \in C^\infty(\Ri^d)$ and $x \in \Ri^d$.
Define $\Phi \colon \Ri^d \times \Ri \to \Ri^d$ by
\[
\Phi(x,t)
= \exp( \sum_{\alpha \in J_s^+(N)} \delta^{-|\alpha|} b(\alpha) \, t^{|\alpha|} X_{[\alpha]})
        \exp(-a_1 t X_{i_1}) \ldots \exp(-a_n t X_{i_n})(x)
\;\;\; .  \]
Then it follows from Lemma~\ref{lhor3.1}.\ref{lhor3.1-3} as in the proof of 
Proposition~\ref{phor301} that there exists an $M' > 0$, depending 
only on the $X_{[\alpha]}$ with $\alpha \in J_s^+(N)$ and $n$, such that 
$|\Phi(x,t) - x| \leq M' \, |t|^{s+1}$ uniformly for all 
$t \in [-1,1]$ and $x \in \Ri^d$.
Replacing $t$ by $\delta \, t$ and $a_i$ by $a_i \delta^{-1}$ yields the proposition.\hfill$\Box$

\section{Subellipticity estimates}\label{Shorf3}

In   this section we prove Theorem~\ref{thor1}.
The proof  follows closely  H\"ormander's reasoning and  
the  subsequent  discussion should be read in conjunction with Section~4 of \cite{Hor1}.
Throughout the section  $X_1,\ldots, X_N$ are  $C_b^\infty$-vector fields
but  we do not require that they satisfy the uniform
H\"ormander condition  until  Proposition~\ref{phor10}.
Set $H_0 = \sum_{i=1}^N X_i^* \, X_i$ with domain $D(H_0) = W^{\infty,2}(\Ri^d)$.
Then $H_0$ is essentially self-adjoint by Proposition~\ref{pcom15}.
Let $H$ denote the closure of $H_0$.
Alternatively, define the quadratic form $h$ on $L_2(\Ri^d)$ by 
\[
h(\varphi)=\sum^N_{i=1} \|X_i\varphi\|_2^2
\]
and domain  $D(h) = \bigcap^N_{i=1} D(X_i)$.
Then the form $h$ is closed
Let $\widetilde H$ be the positive self-adjoint operator associated with the 
closed quadratic form~$h$.
Then obviously $H_0 \subseteq \widetilde H$ and by uniqueness of self-adjoint
extensions one has $H = \widetilde H$.

H\"ormander's proofs are based on the extensive use of H\"older norms.
Therefore we associate  with each  $C_b^\infty$-vector field $X$  a family of such norms.
Specifically for all $\gamma \in \langle0,1]$ we define the H\"older norm
$\|\cdot\|_{2;X,\gamma}$ by
\[
\|\varphi\|_{2;X,\gamma} 
= \|\varphi\|_2 + \sup_{0 < |t| \leq 1} |t|^{-\gamma} \|e^{tX} \varphi - \varphi\|_2
\]
for all $\varphi \in C_c^\infty(\Ri^d)$.
In addition we introduce the universal H\"older norms
\[
\|\varphi\|_{2;\gamma}
= \|\varphi\|_2 + \sup_{0 < |x| \leq 1} |x|^{-\gamma} \|L(x) \varphi - \varphi\|_2
\]
for all $\varphi \in C_c^\infty(\Ri^d)$, where  $L(\,\cdot\,)$  denotes the left regular representation
of $\Ri^d$ on $L_2(\Ri^d)$.

The next  two lemmas are similar to Lemmas~4.1 and 4.2 in \cite{Hor1}.

\begin{lemma} \label{lhor4}
Let $X$ be a $C_b^\infty$-vector field.
Further let $\psi \in C_b^\infty(\Ri^d)$ and $\gamma \in \langle0,1]$.
Then there exists a $c > 0$ such that 
\[
\|\varphi\|_{2;\psi X,\gamma}
\leq c \, \|\varphi\|_{2;X,\gamma}
\]
for all $\varphi \in C_c^\infty(\Ri^d)$.
\end{lemma}
\proof\
Following H\"ormander's proof of Lemma~4.1 we define 
$\tau \colon \Ri^d \times \Ri \to \Ri$ to be 
the solution for each  $x \in \Ri^d$ of the initial value problem 
\[
{\partial_t}\tau(x,t) = \psi(\exp( \tau(x,t) X)(x))
\;\;\;\; \mbox{and} \;\;\;\;
\tau(x,0) = 0
\;\;\; .  \]
Then it follows from Gronwall's lemma that there are $M,\omega > 0$ such that 
\begin{equation}
|\partial_k \tau(x,t)|
\leq M \, e^{\omega |t|}
\label{elhor4;1}
\end{equation}
for all $k \in \{ 1,\ldots,d \} $, $x \in \Ri^d$ and $t \in \Ri$.
Consequently one calculates as in \cite{Hor1} that 
\begin{eqnarray*}
\|e^{t \psi X} \varphi - \varphi\|_2
& \leq & \\
&&\hspace{-1.5cm}{}2 \,|t|^{-1} \int dx \int_{ \{ \sigma : |\sigma| \leq |t| \} } d\sigma \,
       |\varphi( \exp(\tau(x,t) X)(x)) - \varphi(\exp(\sigma X)(x))|^2
   + 2 \,|t|^{2 \gamma} \, \|\varphi\|_{2;X,\gamma}^2
\end{eqnarray*}
for all $t \in [-1,1] \backslash \{ 0 \} $ and $\varphi \in C_c^\infty(\Ri^d)$.

Fix $t \in [-1,1] \backslash \{ 0 \} $.
Introduce new variables $y = \exp(\sigma X)(x)$ and $w = \tau(x,t) - \sigma$.
Then the Jacobian of the coordinate transformation is given by
\[
J_{x,\sigma} 
= \left| \begin{array}{cc}
     \partial_x \exp(\sigma X)(x) & \partial_\sigma \exp(\sigma X)(x)   \\
     (\partial_x \tau)(x,t)       &  -1
  \end{array} \right|
\;\;\; .  
\]
Since $|\sigma| \leq |t| \leq 1$ it follows from (\ref{elhor4;1}) and 
Lemma~\ref{lhor3.1}.\ref{lhor3.1-3} that there exists an $M > 0$ such that 
$|J_{x,\sigma}| \leq M$ uniformly for all $x \in \Ri^d$ and $\sigma \in [-1,1]$.
Moreover, $|\tau(x,t)| \leq \|\psi\|_\infty \, |t|$, 
so $|w| \leq (1 + \|\psi\|_\infty) \, |t|$.
Hence 
\begin{eqnarray*}
\lefteqn{
|t|^{-1} \int dx \int_{ \{ \sigma : |\sigma| \leq |t| \} } d\sigma \,
       |\varphi( \exp(\tau(x,t) X)(x)) - \varphi(\exp(\sigma X)(x))|^2
} \hspace{10mm} \\*
& \leq & M \, |t|^{-1} \int dy \int_{ \{ w : |w| \leq (1 + \|\psi\|_\infty) \, |t| \} } dw \,
       |\varphi( \exp(w X)(y)) - \varphi(y)|^2  
\leq  M' \, |t|^{2\gamma} \, \|\varphi\|_{2;X,\gamma}^2
\end{eqnarray*}
for all $\varphi \in C_c^\infty(\Ri^d)$, where we used Lemma~\ref{lhor3.1}.\ref{lhor3.1-4}.
The statement of the lemma follows immediately.\hfill$\Box$

\begin{lemma} \label{lhor5}
Let $\Phi \colon \Ri^d \times \langle-2,2\rangle \to \Ri^d$ be a $C^\infty$-function,
$N' \in \Ni$ and $\gamma \in \langle0,1]$.
Suppose there exists an $M > 0$ such that 
\[
|\Phi(x,t) - x| \leq M \, t^{N'}
\;\;\;\; \mbox{and} \;\;\;\;
|\partial_k \Phi(x,t)| \leq M
\]
uniformly for all $x \in \Ri^d$, $t \in [-1,1]$ and $k \in \{ 1,\ldots,d \} $.
Then there exists a $c > 0$ such that 
\[
\int_{\Ri^d} dx \, |\varphi(\Phi(x,t)) - \varphi(x)|^2
\leq c \, |t|^{2N'\gamma} \|\varphi\|_{2;\gamma}^2
\]
for all $\varphi \in C_c^\infty(\Ri^d)$ and $t \in [-1,1]$.
\end{lemma}
\proof\
The proof is similar to  the proof of Lemma~4.2 in \cite{Hor1}, with the same 
modifications as in the proof of Lemma~\ref{lhor4}.\hfill$\Box$

\ruimte

The conclusion of Lemma~\ref{lhor5} can be immediately translated into a bound
on the H\"older norm.

\begin{cor}\label{chor5.1}
Let $X$ be a $C_b^\infty$-vector field
and let $\gamma\in\langle0,1]$.
Then there exists a $c>0$ such that $\|\varphi\|_{2;X,\gamma}\leq c\,\|\varphi\|_{2;\gamma}$
for all $\varphi\in C_c^\infty(\Ri^d)$.
\end{cor}
\proof\
It follows from the Duhamel formula that
\[
|\pi_k(\exp(tX))(x)-\pi_k(x)|\leq \int^t_0ds\,|(X\pi_k)(\exp(sX)(x)|\leq \|X\pi_k\|_\infty\,|t|
\]
for all $k\in\{1,\ldots,d\}$.
Therefore
\[
|(\exp(tX)(x)-x|\leq |t| \sum_{k=1}^d \|X \pi_k\|_\infty
\]
for all $t\in\Ri$ and $x\in\Ri^d$.
Then the corollary follows from Lemma~\ref{lhor3.1}.\ref{lhor3.1-3} and Lemma~\ref{lhor5}
applied with $\Phi(x,t)=\exp(tX)(x)$ and $N'=1$.\hfill$\Box$

\ruimte

Proposition~\ref{phor301} immediately yields a global version of H\"ormander's Lemma~4.5.

\begin{lemma} \label{lhor7}
Let $Y_1$ and $Y_2$ be $C_b^\infty$-vector fields,
$\gamma \in \langle0,1]$ and $N \in \Ni \backslash \{ 1 \} $.
Let $Z_2,\ldots,Z_N$ be as in Proposition~{\rm \ref{phor301}}.
Then there exists a $c > 0$ such that 
\[
\|e^{t(Y_1 + Y_2)} \varphi - \varphi\|_2
\le c \,\Big( \|e^{t Y_1} \varphi - \varphi\|_2 + \|e^{t Y_2} \varphi - \varphi\|_2
             + \sum_{j=2}^N \|e^{t^j Z_j} \varphi - \varphi\|_2
             + |t|^{\gamma(N+1)} \|\varphi\|_{2;\gamma} \Big)
\]
for all $\varphi \in C_c^\infty(\Ri^d)$ and $t \in [-1,1]$.
\end{lemma}
\proof\
Define $\Phi \colon \Ri^d \times \Ri \to \Ri^d$ by
\[
\Phi(x,t) 
= \exp( t(Y_1 + Y_2) ) \exp(-t Y_1) \exp(-t Y_2) \exp(-t^2 Z_2) \ldots
            \exp(- t^N Z_N)(x) 
\;\;\; .  \]
If $c > 0$ is as in Proposition~\ref{phor301} then 
it follows  that
$|\Phi(x,t) - x| \leq c \, |t|^{N+1}$
for all $x \in \Ri^d$ and $t \in [-1,1]$.
Secondly, 
\[
\sup_{k \in \{ 1,\ldots,d \} }
\sup_{x \in \Ri^d}
\sup_{t \in [-1,1]}
       |(\partial_k \Phi)(x,t)|
   < \infty
\]
again by Lemma~\ref{lhor3.1}.\ref{lhor3.1-3}.
Hence by Lemma~\ref{lhor5} it follows that there is a $c_1 > 0$ such that 
\[
\int dx \, |\varphi(\Phi(x,t)) - \varphi(x)|^2
\leq c_1^2 \, |t|^{2(N+1)\gamma} \|\varphi\|_{2;\gamma}^2
\]
for all $\varphi \in C_c^\infty(\Ri^d)$ and $t \in [-1,1]$.

Next, for all $t \in \Ri$ define $H_t \colon C^\infty(\Ri^d) \to C^\infty(\Ri^d)$ by
$(H_t \varphi)(x) = \varphi(\Phi(x,t))$.
Then $\|H_t \varphi - \varphi\|_2 \leq c_1 \, t^{(N+1) \gamma} \|\varphi\|_{2;\gamma}$
for all $\varphi \in C_c^\infty(\Ri^d)$ and $t \in [-1,1]$.
But 
\[
H_t = e^{-t^N Z_N} \ldots e^{-t^2 Z_2} e^{-t Y_2} e^{-t Y_1} e^{ t(Y_1 + Y_2) }
\]
and 
\[
e^{ t(Y_1 + Y_2) }
= e^{t Y_1} e^{t Y_2} e^{t^2 Z_2} \ldots e^{ t^N Z_N} H_t
\]
for all $t \in \Ri$.
Then the lemma follows from a concertina formula, and
Lemma~\ref{lhor3.1}.\ref{lhor3.1-4}.~\hfill$\Box$

\ruimte

We emphasize that in the next two lemmas it is not necessary for the vector fields
to satisfy the uniform H\"ormander condition.

\begin{lemma} \label{lhor8}
Let $X_1,\ldots,X_{N}$ be $C_b^\infty$-vector fields
and let $\gamma,\delta \in \langle0,1]$.
Then for all $\alpha\in J^+(N)$ there exists $c_1,c_2 > 0$ such that 
\[
\|\varphi\|_{2;X_{[\alpha]},\gamma |\alpha|^{-1}}
\leq c_1 \sum_{j=1}^{N} \|\varphi\|_{2;X_j,\gamma} 
   + c_2 \, \|\varphi\|_{2;\delta}
\]
for all $\varphi \in C_c^\infty(\Ri^d)$.
\end{lemma}
\proof\
Let $\alpha \in J^+(N)$.
Fix $s \in \Ni$ with $s \geq |\alpha| \vee \gamma \, \delta^{-1}$.
By Proposition~\ref{phor843} applied with $\delta = 1$ there exist
$n \in \{ 1,\ldots,3(2d)^s \} $, $i_1,\ldots,i_n \in \{ 1,\ldots,N \} $,
$M,M' > 0$ and $a_1,\ldots,a_n \in [-M,M]$ such that if 
$\Phi \colon \Ri^d \times \Ri \to \Ri^d$ is given by
\[
\Phi(x,t)
= \exp( t^{|\alpha|} X_{[\alpha]})
        \exp(-a_1 t X_{i_1}) \ldots \exp(-a_n t X_{i_n})(x)
\]
then $|\Phi(x,t) - x| \leq M' \, |t|^{s+1}$ uniformly for all 
$t \in [-1,1]$ and $x \in \Ri^d$.
Moreover, by Lemma~\ref{lhor3.1}.\ref{lhor3.1-3}, there exists an $M''>0$ such that
$|\partial_k \Phi(x,t)| \leq M''$
uniformly for all $x \in \Ri^d$, $t \in [-1,1]$ and $k \in \{ 1,\ldots,d \} $.
For all $t \in \Ri$ define $H_t \colon C^\infty(\Ri^d) \to C^\infty(\Ri^d)$ by
$(H_t \varphi)(x) = \varphi(\Phi(x,t))$.
Then by Lemma~\ref{lhor5}  there is a $c > 0$ such that 
\[
\|H_t \varphi - \varphi\|_2
\leq c \, |t|^{(s+1) \delta} \|\varphi\|_{2;\delta} 
\leq c \, |t|^\gamma \|\varphi\|_{2;\delta }
\]
for all $t \in [-1,1]$ and $\varphi \in C_c^\infty(\Ri^d)$.
But 
\[
e^{t^{|\alpha|} X_{[\alpha]}} 
= e^{a_1 t X_{j_1}} \ldots e^{a_n t X_{i_n}} H_t
\;\;\; .  \]
Therefore Lemma~\ref{lhor3.1}.\ref{lhor3.1-4} implies that there is a $c' > 0$ such that 
\begin{eqnarray*}
\|e^{t^{|\alpha|} X_{[\alpha]}} \varphi - \varphi\|_2
& \leq & c' \,\Big( \sum_{l=1}^n \|e^{a_l t X_{i_l}} \varphi - \varphi\|_2
             + \|H_t \varphi - \varphi\|_2 
        \Big)  \\
& \leq & c' \,\Big( |t|^\gamma \sum_{l=1}^n \|\varphi\|_{2;a_l X_{i_l},\gamma}
             + c \, |t|^\gamma \|\varphi\|_{2;\delta }
        \Big)
\end{eqnarray*}
for all $t \in [-1,1]$ and $\varphi \in C_c^\infty(\Ri^d)$.
Then the lemma follows from Lemma~\ref{lhor4}.\hfill$\Box$

\begin{lemma} \label{lhor9}
Let $X_1,\ldots,X_{N}$ be $C_b^\infty$-vector fields.
For all $k \in \Ni$ set
\[
{\cal D}^{(k)}
= \spann \{ \psi X_{[\alpha]} : \psi \in C_b^\infty(\Ri^d) , \; \alpha \in J_k^+(N)\} 
\;\;\; .  \]
Then for all $\delta,\gamma \in \langle0,1]$, $k \in \Ni$ and $X \in {\cal D}^{(k)}$
there exists a $c > 0$ such that 
\begin{equation}
\|\varphi\|_{2;X,\gamma k^{-1}} 
\leq c \Big( \sum_{j=1}^{N} \|\varphi\|_{2;X_j,\gamma} 
             + \|\varphi\|_{2;\delta}
       \Big)
\label{elhor9;1}
\end{equation}
for all $\varphi \in C_c^\infty(\Ri^d)$.
\end{lemma}
\proof\
Fix $\delta \in \langle0,1]$.
If $\gamma \leq \delta k$ then (\ref{elhor9;1}) follows from Corollary~\ref{chor5.1}.
For all $n \in \Ni$ let $P(n)$ be the following hypothesis.
\begin{quote}
For all $k \in \Ni$, $\gamma \in \langle0,\delta k 2^{n-2} \wedge 1]$ and 
$X \in {\cal D}^{(k)}$ 
there exists a $c > 0$ such that 
\[
\|\varphi\|_{2;X,\gamma k^{-1}} 
\leq c\, \Big( \sum_{j=1}^{N} \|\varphi\|_{2;X_j,\gamma} 
             + \|\varphi\|_{2;\delta}
       \Big)
\]
for all $\varphi \in C_c^\infty(\Ri^d)$.
\end{quote}
Then $P(1)$ is valid.
Let $n \in \Ni$ and suppose that $P(n)$ is valid.
Let $k \in \Ni$ and $\gamma \in \langle0,\delta k 2^{n-1} \wedge 1]$.
Consider
\begin{eqnarray*}
V & = & \{ X \in {\cal D}^{(k)} : \mbox{there exists a } c > 0 \mbox{ such that }
   \\*
& & \hspace{30mm}
       \|\varphi\|_{X,\gamma k^{-1}} 
           \leq c \,\Big( \sum_{j=1}^{N} \|\varphi\|_{2;X_j,\gamma} 
             + \|\varphi\|_{2;\delta} \Big)
        \mbox{ for all } \varphi \in C_c^\infty(\Ri^d) 
    \} 
\;\;\; .  
\end{eqnarray*}
If $\alpha \in J_k^+(N)$ then $X_{[\alpha]} \in V$ by Lemma~\ref{lhor8}.
Moreover, if in addition $\psi \in C_b^\infty(\Ri^d)$ then 
$\psi X_{[\alpha]} \in V$ by Lemma~\ref{lhor4}.
So it remains to show that $V$ is a vector space.
But that follows from Lemma~\ref{lhor7} and the induction hypothesis.\hfill$\Box$

\ruimte

The next proposition is the first application of the uniform
H\"ormander condition.

\begin{prop} \label{phor10}
Let $X_1,\ldots,X_{N}$ be $C_b^\infty$-vector fields satisfying the  uniform
H\"ormander condition of order $r$ on $\Ri^d$.
Then for all $\gamma \in \langle0,1]$ there exists a $c > 0$ such that 
\[
\|\varphi\|_{2;\gamma r^{-1}} 
\leq c \,\Big( \sum_{j=1}^{N} \|\varphi\|_{2;X_j,\gamma} + \|\varphi\|_2 \Big)
\]
for all $\varphi \in C_c^\infty(\Ri^d)$.
\end{prop}
\proof\
It follows from Lemma~\ref{lhor9} that for all $X \in {\cal D}^{(r)}$ there exists a $c > 0$ 
such that 
\[
\|\varphi\|_{2;\gamma r^{-1}} 
\leq c\, \Big( \sum_{j=1}^{N} \|\varphi\|_{2;X_j,\gamma} + \|\varphi\|_{2;\gamma (2r)^{-1}} \Big)
\]
for all $\varphi \in C_c^\infty(\Ri^d)$,
where ${\cal D}^{(r)}$ is as in Lemma~\ref{lhor9}.
By the uniform H\"ormander condition one has 
$\partial_i \in {\cal D}^{(r)}$ for all $i \in \{ 1,\ldots,d \} $.
Hence there is a $c > 0$ such that 
\[
\|\varphi\|_{2;\gamma r^{-1}} 
\leq d \sum_{i=1}^d \|\varphi\|_{2;\partial_i,\gamma r^{-1}} 
\leq c \,\Big( \sum_{j=1}^{N} \|\varphi\|_{2;X_j,\gamma} + \|\varphi\|_{2;\gamma (2r)^{-1}} \Big)
\]
for all $\varphi \in C_c^\infty(\Ri^d)$.
But there is a $c_1 > 0$ such that 
\[
\|\varphi\|_{2;\gamma (2r)^{-1}}
\leq \varepsilon \,\|\varphi\|_{2;\gamma r^{-1}} + c_1 \, \varepsilon^{-1} \|\varphi\|_2 
\]
for all $\varphi \in C_c^\infty(\Ri^d)$ and $\varepsilon > 0$.
Choosing $\varepsilon = (2c)^{-1}$ one deduces that 
\[
\|\varphi\|_{2;\gamma r^{-1}} 
\leq 2c\, \Big( \sum_{j=1}^{N} \|\varphi\|_{2;X_j,\gamma} + 2 c \, c_1 \, \|\varphi\|_2 \Big)
\]
for all $\varphi \in C_c^\infty(\Ri^d)$.\hfill$\Box$

\ruimte

For the last part of the proof of Theorem~\ref{thor1} we need some additional  interpolation
spaces.
The proof relies on an extrapolation, interpolation and a similar
extrapolation argument.
If $r = 1$ then $H_0$ is strongly elliptic and the theorem is well known. 
So we may assume that $r \geq 2$.

If $L$ is the generator of a continuous semigroup $S$ on $L_2(\Ri^d)$,
$p \in [1,\infty]$ and $\gamma \in \langle0,1]$ define the functions
$\|\cdot\|_{\gamma,p,S}, \;\; 
\|\cdot\|_{\gamma,p,S}', \;\;
\|\cdot\|_{\gamma,p,L} \colon L_2(\Ri^d) \to [0,\infty]$
by
\begin{eqnarray*}
\|\varphi\|_{\gamma,p,S}
& = & \|\varphi\|_2 
   + \Big( \int_0^1 dt\, t^{-1}\,\Big|\, t^{-\gamma} \, \|(I - S_t) \varphi\|_2\, \Big|^p \Big)^{1/p} \;\;\;, \\
\|\varphi\|_{\gamma,p,S}'
& = & \|\varphi\|_2 
   + \Big( \int_0^1dt\, t^{-1}\,\Big|\, t^{-\gamma} \, \|(I - S_t)^2 \varphi\|_2 \,\Big|^p \Big)^{1/p} \;\;\;, \\
\|\varphi\|_{\gamma,p,L}
& = & \|\varphi\|_2 
   + \Big( \int_0^1 dt\, t^{-1}\, \Big| \,t^{-\gamma} \, \kappa_t(\varphi) \,\Big|^p \Big)^{1/p}  \;\;\;,
\end{eqnarray*}
if $p < \infty$, where 
\[
\kappa_t(\varphi)
= \inf \{ \|\varphi - \varphi_1\|_2 + t \, \|L \varphi_1\|_2 : \varphi_1 \in D(L) \} 
\]
 with obvious modifications if $p = \infty$.
Define the interpolation spaces 
\begin{eqnarray*}
\cx_{\gamma,p,S}
& = & \{ \varphi \in L_2(\Ri^d) : \|\varphi\|_{\gamma,p,S} < \infty \} \\
\cx_{\gamma,p,S}'
& = & \{ \varphi \in L_2(\Ri^d) : \|\varphi\|_{\gamma,p,S}' < \infty \} 
\end{eqnarray*}
with norms $\|\cdot\|_{\gamma,p,S}$ and $\|\cdot\|_{\gamma,p,S}'$.
If $\cx$ and $\cy$ are two Banach spaces which are embedded in a locally convex Hausdorff space
denote by $(\cx,\cy)_{\gamma,p,K}$ the interpolation space with respect to the K-method.
Then 
\[
(L_2(\Ri^d),D(L))_{\gamma,p,K} 
= \{ \varphi \in L_2(\Ri^d) : \|\varphi\|_{\gamma,p,L} < \infty \} 
\]
and the norm is equivalent to $\|\cdot\|_{\gamma,p,L}$.
If $S$ is a continuous semigroup then it follows from \cite{BB}, Theorem~3.4.2 and Corollary~3.4.9,
that the spaces $\cx_{\gamma,p,S}$, $\cx_{\gamma,p,S}'$ and $(L_2(\Ri^d),D(L))_{\gamma,p,K}$ 
are equal with equivalent norms if $\gamma < 1$.
Moreover, if $S$ is merely continuous, $p = \infty$ and $\gamma = 1$ then 
$D(L) \subset \cx_{1,\infty,S}$ and the embedding is continuous.
If $L$ is a positive self-adjoint operator and $p = 2$ then a much better result
is valid:
\[
D(L^\gamma) = \cx_{\gamma,2,S}
\]
and the norms are equivalent (see \cite{ER1}, Lemma~7.1).

As in Section~\ref{Shor3a}
we set $L = I +\Delta$ and let $S$ be the semigroup
generated by~$L$.

\begin{lemma} \label{lhor337}
$D(H) \subset (L_2(\Ri^d),D(L))_{r^{-1},\infty,K}$ and the embedding is continuous.
\end{lemma}
\proof\
It follows from \cite{ER1}, Theorem~3.2,
that the norms $\|\cdot\|_{2,\delta}$ and $\cx_{2^{-1} \delta,\infty,S}$
are equivalent for all $\delta \in \langle0,1\rangle$.
Moreover, $D(X_i) \subset \cx_{1,\infty,X_i}$ and the embedding is continuous.
Hence it follows from Proposition~\ref{phor10}, applied with $\gamma = 1$,
that there is a $c_1 > 0$ such that 
\begin{equation}
\sup_{0 < t \leq 1} t^{-1/r} \|(I - S_t) \varphi\|_2^2
\leq c_1 \Big( \|\varphi\|_2^2 + (\varphi,H\varphi) \Big)
\label{elhor337;1}
\end{equation}
for all $\varphi \in C_c^\infty(\Ri^d)$.
Then by density (\ref{elhor337;1}) is valid for all $\varphi \in W^{2,2}(\Ri^d)$.
Next let $c > 0$ be as in Corollary~\ref{chord2.1}.\ref{chord2.1-2}.
Set $\tau = (2r)^{-1}$ and let $t \in \langle0,1]$.
Choosing $A = c_1 \, (I+H_0)$ and $B = t^{-\tau} (I - S_t)$ it follows from (\ref{elhor337;1})  that 
$(\varphi, A \varphi) \geq \|B \varphi\|^2$ for all $\varphi \in W^{2,2}(\Ri^d)$.
Moreover, 
\begin{eqnarray*}
|(\varphi, [B, [B,A]] \varphi)|
= c_1 \, t^{-2 \tau} \, |(\varphi, [S_t, [S_t,H_0]] \varphi)|  
\leq c \, c_1 \, \|B\varphi\|_2^2
\leq 2^{-1} \|B^2 \varphi\|_2^2 + c^2 \, c_1^2 \, \|\varphi\|_2^2
\end{eqnarray*}
for all $\varphi \in W^{\infty,2}(\Ri^d)$ by Corollary~\ref{chord2.1}.
Therefore the assumptions of Lemma~\ref{lcom2}
are valid with $\cd=W^{\infty,2}(\Ri^d)$ uniformly for all
 $t \in \langle0,1]$.
Hence 
\begin{eqnarray*}
c \, \|(I + H_0) \varphi\|_2
 =  \|A \varphi\|_2
 \geq  2^{-1} \|B^2 \varphi\|_2 - c \, c_1 \,\|\varphi\|_2  
 =  2^{-1} t^{-2\tau} \|(I - S_t)^2 \varphi\|_2 - c \, c_1 \, \|\varphi\|_2 
\end{eqnarray*}
uniformly for all $\varphi \in W^{\infty,2}(\Ri^d)$ and $t \in \langle0,1]$.
Therefore 
\[
c \, \|(I + H_0) \varphi\|_2
\geq 2^{-1} \|\varphi\|_{2\tau,\infty,S}' - (c \, c_1 + 1) \, \|\varphi\|_2
\]
for all $\varphi \in W^{\infty,2}(\Ri^d)$.
Since $W^{\infty,2}(\Ri^d)$ is a core of $H_0$ and $\cx'_{2\tau,\infty,S}$ is complete it follows that 
$D(H) \subset \cx'_{2\tau,\infty,S}$ and the embedding is continuous. 
Finally, the lemma follows because 
$\cx_{2\tau,\infty,S}' = \cx_{2\tau,\infty,L}$, with equivalent norms.\hfill$\Box$

\ruimte

It follows from Lemma~\ref{lhor337} that 
\[
D(H) \subset (L_2(\Ri^d),D(L))_{r^{-1},\infty,K} 
\]
and the embedding is continuous.
Hence by interpolation
\[
(L_2(\Ri^d),D(H))_{2^{-1},2,K} 
\subset (L_2(\Ri^d), (L_2(\Ri^d),D(L))_{r^{-1},\infty,K} )_{2^{-1},2,K} 
\;\;\; .  \]
But by the reiteration theorem, \cite{BB} Theorem~3.2.20, one has
\[
(L_2(\Ri^d), (L_2(\Ri^d),D(L))_{r^{-1},\infty,K} )_{2^{-1},2,K} 
= (L_2(\Ri^d),D(L))_{\tau,2,K} 
\]
with equivalent norms, where again $\tau = (2r)^{-1}$.
Moreover, 
\[
(L_2(\Ri^d),D(L))_{\tau,2,K} = D(L^\tau)
\;\;\;\;
\mbox{and} 
\;\;\;\;
(L_2(\Ri^d),D(H))_{2^{-1},2,K} = D(H^{1/2})
\]
with equivalent norms.
So $D(H^{1/2}) \subset D(L^\tau)$ and the embedding is continuous.
Therefore there exists a $c > 0$ such that 
\[
\|\Delta^\tau\varphi\|_2^2\leq \|L^\tau \varphi\|_2^2
\leq c \,\Big( \|\varphi\|_2^2 + \|H^{1/2} \varphi\|_2^2 \Big)
= c \,\Big( \|\varphi\|_2^2 + (\varphi, H\varphi) \Big)
\]
for all $\varphi \in D(H)$.
Then Theorem~\ref{thor1} is a corollary of Theorem~\ref{tcom31}.\hfill$\Box$

\section{The  Uniform H\"ormander Condition}\label{Shorf2}

Let $X_1,\ldots,X_N$ be $C_b^\infty$-vector fields on $\Ri^d$.
We conclude by deriving several characterizations of the uniform version of the H\"ormander condition.

Each vector field $X_i$ can be expressed as a partial differential operator
$X_i=\sum^d_{k=1} a_{ik}\,\partial_k$ 
with coefficients $a_{ik} = X_i \pi_k\in C_b^\infty(\Ri^d)$, where $\pi_k$ denotes the projection on the 
$k$-th coordinate.
The multi-commutator
$X_{[\alpha]}$ is also a $C_b^\infty$-vector field with coefficients 
$a_{\alpha\,k} = X_{[\alpha]} \pi_k \in C_b^\infty(\Ri^d)$.
Explicitly 
$X_{[\alpha]}=\sum^d_{k=1} a_{\alpha \,k}\,\partial_k$.
Then for all $r\in\Ni$ and  all $i,j\in\{1,\ldots,d\}$ define
\[
c^{(r)}_{ij}
=\sum_{\alpha \in J_r^+(N)}
      a_{\alpha \,i}\,a_{\alpha \,j}
\]
and set $C^{(r)}=(c^{(r)}_{ij})$.
The matrix $C^{(r)}$ is real symmetric and positive semidefinite.
In particular the operator $ H_0$ given by (\ref{etcom11;2}) is a second-order 
operator in divergence form with the matrix of coefficients
$C^{(1)}=(\sum^d_{k=1}a_{ki} \, a_{kj})$.

\begin{prop}\label{phor2}
Let $X_1,\ldots,X_N$ be $C_b^\infty$-vector fields.
For all $\alpha \in J^+(N)$ and $x \in \Ri^d$ let $a_\alpha(x) \in \Ri^d$ be such that 
$X_{[\alpha]}=\sum^d_{k=1} a_{\alpha \,k}\,\partial_k$.
Moreover, fix $r \in \Ni$.
Then the following statements are equivalent.
\begin{tabel}
\item \label{phor2-2}
The vector fields $X_1,\ldots,X_N$   satisfy the uniform H\"ormander
condition  of order $r$ on~$\Ri^d$.
\item \label{phor2-1}
There exists a $\sigma>0$ such that $C^{(r)}(x)\geq\sigma I$ uniformly for all $x\in\Ri^d$.
\item \label{phor2-2.1}
There exists an $M > 0$ such that for all $x \in \Ri^d$,
$i \in \{ 1,\ldots,d \} $ and $\alpha \in J_r^+(d')$ there exists a 
$\lambda_\alpha \in [-M,M]$ such that 
\[
e_i = \sum_{\alpha \in J_r^+(d')} \lambda_\alpha \, a_\alpha(x)
\;\;\; ,  \]
where $e_i$ is the unit vector in the $i$-th direction.
\item \label{phor2-3}
There exists a $\sigma > 0$ such that 
\[
\Vol \Big\{ \sum_{\alpha \in J^+_r(N)} \lambda_\alpha \, a_\alpha(x)
          : |\lambda_\alpha| \leq 1 \mbox{ for all } \alpha \in J^+_r(N) \Big\} 
\geq \sigma
\]
uniformly for all $x \in \Ri^d$.
\item \label{phor2-4}
There exists a $\sigma > 0$ such that for all $x \in \Ri^d$ there are multi-indices
$\alpha_1,\ldots,\alpha_d \in J^+_r(N)$ such that 
\[
|\det( (X_{[\alpha_i]} \pi_j)(x) )|
= |\det(a_{\alpha_1}(x),\ldots,a_{\alpha_d}(x))| 
\geq \sigma
\;\;\; .  
\]
\end{tabel}
\end{prop}
\proof\
\ref{phor2-2}$\Rightarrow$\ref{phor2-2.1}.
It follows from Statement~\ref{phor2-2} that for all $i \in \{ 1,\ldots,d \} $  
and $\alpha \in J^+_r(N)$ there
are $\psi_{i\, \alpha} \in C_b^\infty(\Ri^d)$
such that 
\[
\partial_i = \sum_{\alpha\in J^+_r(N)} \psi_{i\, \alpha} X_{[\alpha]}
\]
for all $i \in \{ 1,\ldots,d \} $.
Then $e_i = \sum_{\alpha\in J^+_r(N)} \psi_{i\, \alpha}(x) \, a_\alpha(x)$
for all $x \in \Ri^d$
and Statement~\ref{phor2-2.1} follows with
$M = \max_{i \in \{ 1,\ldots,d \} } \max_{\alpha\in J^+_r(N)} \|\psi_{i\, \alpha}\|_\infty$.

\smallskip

\ref{phor2-2.1}$\Rightarrow$\ref{phor2-3}.
Let $M > 0$ be as in Statement~\ref{phor2-2.1}.
Then 
\[
\Big\{ \sum_{i=1}^d \lambda_i \, e_i : 0 \leq \lambda_i \leq (dM)^{-1} \mbox{ for all } i \Big\} 
\subseteq \Big\{ \sum_{\alpha\in J^+_r(N)} \lambda_\alpha \, a_\alpha(x)
                 : |\lambda_\alpha| \leq 1 \mbox{ for all } \alpha \Big\} 
\]
for all $x \in \Ri^d$.
Therefore 
\begin{eqnarray*}
\Vol  \Big\{ \sum_{\alpha\in J^+_r(N)} \lambda_\alpha \, a_\alpha(x)
                 : |\lambda_\alpha| \leq 1 \mbox{ for all } \alpha \Big\} 
& \geq & \Vol \Big\{ \sum_{i=1}^d \lambda_i \, e_i : 0 \leq \lambda_i \leq (dM)^{-1} 
\mbox{ for all } i \Big\}  \\
& = & (dM)^{-d}
\end{eqnarray*}
for all $x \in \Ri^d$ and Statement \ref{phor2-3} follows.

\smallskip

\ref{phor2-3}$\Rightarrow$\ref{phor2-4}.
Fix $x \in \Ri^d$.
By Lemma~3.1.2 in \cite{Smu1} there are $\alpha_1,\ldots,\alpha_d \in J_r^+(N)$ 
and for all $\alpha \in J_r^+(N)$ and $k \in \{ 1,\ldots,d \} $ there are
$\lambda_{\alpha \,k} \in \Ri$ with $|\lambda_{\alpha\, k}| \leq 2^{L-d}$
such that 
\[
a_\alpha(x) = \sum_{k=1}^d \lambda_{\alpha \,k} \, a_{\alpha_k}(x)
\]
 where $L = \card J_r^+(N)$.
Then 
\[
\Big\{ \sum_{\alpha\in J^+_r(N)} \lambda_\alpha \, a_\alpha(x)
                 : |\lambda_\alpha| \leq 1 \mbox{ for all } \alpha \Big\} 
\subseteq 2^{L-d} L \Big\{ \sum_{k=1}^d \lambda_k \, a_{\alpha_k}(x) : |\lambda_k| \leq 1
                                           \mbox{ for all } k \Big\}
\;\;\; .  \]
Therefore 
\begin{eqnarray*}
|\det(a_{\alpha_1}(x),\ldots,a_{\alpha_d}(x))| 
& = & 2^{-d} \Vol  \Big\{ \sum_{k=1}^d \lambda_k \, a_{\alpha_k}(x) : |\lambda_k| \leq 1
                                           \mbox{ for all } k \Big\}  \\
& \geq & 2^{-L} \, L^{-1} \Vol \Big\{ \sum_{\alpha \in J_r^+(N)} \lambda_\alpha \, a_\alpha(x)
                 : |\lambda_\alpha| \leq 1 \mbox{ for all } \alpha \Big\}  
\end{eqnarray*}
Thus Statement~\ref{phor2-3} implies Statement~\ref{phor2-4}.

\smallskip

\ref{phor2-4}$\Rightarrow$\ref{phor2-1}.
Let $\sigma > 0$ be as in Statement~\ref{phor2-4}.
Fix $x \in \Ri^d$.
Then there are $\alpha_1,\ldots,\alpha_d \in J^+_r(N)$
such that 
\[
|\det(a_{\alpha_1}(x),\ldots,a_{\alpha_d}(x))| \geq \sigma
\;\;\; .  \]
For all $k,l \in \{ 1,\ldots,d \} $ set 
$d_{kl} = \sum_{j=1}^d a_{\alpha_j \,k}(x) \, a_{\alpha_j \,l}(x)$ and $D = (d_{kl})$.
Then 
\[
\det D 
= \Big( \det(a_{\alpha_1}(x),\ldots,a_{\alpha_d}(x)) \Big)^2 
\geq \sigma^2
\;\;\; .  \]
Moreover, 
\[
C^{(r)}
\geq D
\geq \|D\|^{-(d-1)} \, (\det D) \, I
\geq \|D\|^{-(d-1)} \, \sigma^2 \, I
\]
where $\|D\|$ is the norm of the matrix $D$.
Since the coefficients $a_\alpha$ are uniformly bounded Statement~\ref{phor2-1} 
follows.

\smallskip

\ref{phor2-1}$\Rightarrow$\ref{phor2-2}.
If $X=\sum^d_{k=1}a_k\,\partial_k$  then since $C^{(r)}$ is invertible one computes that
\[
X=\sum_{\alpha\in J^+_r(N)}((C^{(r)})^{-1}a,a_\alpha)\,
X_{[\alpha]}
\;\;\;.
\]
But the condition $C^{(r)}\geq\sigma I$
implies that the coefficients of the matrix $(C^{(r)})^{-1}$ are in $C_b^\infty(\Ri^d)$.
\hfill$\Box$

\ruimte

Statement~\ref{phor2-1} of 
Proposition~\ref{phor2} is the formulation of the uniform H\"ormander condition used by
Kusuoka and Stroock in Section~3 {\em et seq.} of 
\cite{KuS3} and again in their  analysis of long time behaviour in 
\cite{KuS} (see Theorems~3.20 and 3.24).
The determinant identified in Statement~\ref{phor2-4}
of Proposition~\ref{phor2}
 plays a ubiquitous role 
in the analysis of Nagel, Stein and Wainger \cite{NSW}
and was also identified by Jerison  as an important parameter in the Poincar\'e inequality
(see \cite{Jer}, Condition~(2.3c)  on page~505).

Finally we note that  for operators $H_0$ with $C^\infty$-coefficients Fefferman and Phong have shown that the  subellipticity condition
(\ref{esub1.1}) is locally equivalent to a property of the geometry associated with $H_0$.
 Nagel, Stein and Wainger \cite{NSW} have then analyzed in detail the local geometry for operators
 (\ref{etcom11;2})
 constructed from vector fields satisfying the local H\"ormander condition.
 One could expect that there are global analogues of these results.
 In a separate paper we will indeed extend the conclusions of Nagel, Stein and Wainger 
and  obtain uniform properties of the geometry, properties such as volume doubling,  if the vector fields satisfy the uniform H\"ormander condition.

\subsection*{Acknowledgement}
This work  was  supported by an Australian Research Council (ARC) Discovery Grant DP 0451016.
A part of this work was carried out during a visit of the first named
author to the Australian National University and of the second named author to 
the University of Auckland.


\begin{thebibliography}{CFKS87}

\bibitem[BuB]{BB}
{\sc Butzer, P.L. {\rm and} Berens, H.}, {\em Semi-groups of operators and
  approximation}.
\newblock Die Grundlehren der mathematischen Wissenschaften 145.
  Springer-Verlag, Berlin etc., 1967.

\bibitem[CFKS]{CFKS}
{\sc Cycon, H.L., Froese, R.G., Kirsch, W. {\rm and} Simon, B.}, {\em
  Schr{\"o}dinger operators with application to quantum mechanics and global
  geometry}.
\newblock Springer-Verlag, Berlin etc., 1987.

\bibitem[DrS]{DrS}
{\sc Driessler, W. {\rm and} Summers, S.J.}, On commutators and
  selfadjointness.
\newblock {\em Lett.\ Math.\ Phys.} {\bf 7} (1983),  319--326.

\bibitem[ElR1]{ER1}
{\sc Elst, A.F.M. ter {\rm and} Robinson, D.W.}, Subelliptic operators on
  Lie groups: regularity.
\newblock {\em J. Austr.\ Math.\ Soc.\ {\rm (}Series A{\rm )}} {\bf 57} (1994),
   179--229.

\bibitem[ElR2]{ER19}
\leavevmode\vrule height 2pt depth -1.6pt width 23pt, Second-order strongly
  elliptic operators on Lie groups with H{\"o}lder continuous coefficients.
\newblock {\em J. Austr.\ Math.\ Soc.\ {\rm (}Series A{\rm )}} {\bf 63} (1997),
   297--363.

\bibitem[Far]{Far}
{\sc Faris, W.G.}, {\em Self-adjoint operators}.
\newblock Lect.\ Notes in Math. 433. Springer-Verlag, Berlin etc., 1975.

\bibitem[FeP]{FP}
{\sc Fefferman, C. {\rm and} Phong, D.H.}, Subelliptic eigenvalue problems.
\newblock In {\em Conference on harmonic analysis in honor of Antoni Zygmund},
  Wadsworth Math.\ Ser.,  590--606. Wadsworth, Belmont, CA, 1983.

\bibitem[GlJ1]{GJ1}
{\sc Glimm, J. {\rm and} Jaffe, A.}, The $\lambda\phi_2^4$ quantum field theory
  without cutoffs. IV. Perturbations of the Hamiltonian.
\newblock {\em J. Math.\ Phys.} {\bf 13} (1972),  1568--1584.

\bibitem[GlJ2]{GJ}
\leavevmode\vrule height 2pt depth -1.6pt width 23pt, {\em Quantum Physics. A
  Functional integral point of view}.
\newblock Springer-Verlag, New York etc., 1981.

\bibitem[H{\"o}r]{Hor1}
{\sc H{\"o}rmander, L.}, Hypoelliptic second order differential equations.
\newblock {\em Acta Math.} {\bf 119} (1967),  147--171.

\bibitem[Jer]{Jer}
{\sc Jerison, D.}, The Poincar{\'e} inequality for vector fields satisfying
  H{\"o}rmander's condition.
\newblock {\em Duke Math.\ J.} {\bf 53} (1986),  503--523.

\bibitem[KuS1]{KuS3}
{\sc Kusuoka, S. {\rm and} Stroock, D.}, Applications of the Malliavin
  calculus. III.
\newblock {\em J. Fac.\ Sci.\ Univ.\ Tokyo Sect.\ IA Math.} {\bf 34} (1987),
  391--442.

\bibitem[KuS2]{KuS}
\leavevmode\vrule height 2pt depth -1.6pt width 23pt, Long time estimates for
  the heat kernel associated with a uniformly subelliptic symmetric second
  order operator.
\newblock {\em Ann.\ Math.} {\bf 127} (1988),  165--189.

\bibitem[LuP]{LuP}
{\sc Lumer, G. {\rm and} Phillips, R.S.}, Dissipative operators in a Banach
  space.
\newblock {\em Pacific J. Math.} {\bf 11} (1961),  679--698.

\bibitem[NSW]{NSW}
{\sc Nagel, A., Stein, E.M. {\rm and} Wainger, S.}, Balls and metrics defined
  by vector fields I: basic properties.
\newblock {\em Acta Math.} {\bf 155} (1985),  103--147.

\bibitem[ReS]{RS2}
{\sc Reed, M. {\rm and} Simon, B.}, {\em Methods of modern mathematical physics
  II. Fourier analysis, self-adjoint\-ness}.
\newblock Academic Press, New York etc., 1975.

\bibitem[Rob]{Rob7}
{\sc Robinson, D.W.}, Commutator theory on Hilbert space.
\newblock {\em Can.\ J. Math.} {\bf 34} (1987),  1235--1280.

\bibitem[RoS]{RS}
{\sc Rothschild, L.P. {\rm and} Stein, E.M.}, Hypoelliptic differential
  operators and nilpotent groups.
\newblock {\em Acta Math.} {\bf 137} (1976),  247--320.

\bibitem[Smu]{Smu1}
{\sc Smulders, C.M.P.A.}, {\em Reduced heat kernels on homogeneous spaces}.
\newblock PhD thesis, Eindhoven University of Technology, The Netherlands,
  2000.
\newblock See http://alexandria.tue.nl/extra2/200000153.pdf.

\end{thebibliography}
\end{document}